\pdfoutput=1

\documentclass[final]{siamltex}
\usepackage{eurosym}
\usepackage{algorithm}
\usepackage{algpseudocode}
\usepackage{amsfonts}
\usepackage{amsmath}
\usepackage{amssymb}
\usepackage{enumerate}
\usepackage{graphicx}
\usepackage{multirow}
\usepackage{url}
\usepackage{verbatim}
\usepackage[pdftex]{thumbpdf}
\usepackage[pdftex,
pdfstartview=FitH,bookmarks=true,bookmarksnumbered=true,bookmarksopen=true,hypertexnames=false,breaklinks=true,
colorlinks=true,  linkcolor=blue,anchorcolor=blue,
citecolor=blue,filecolor=blue,  menucolor=blue]{hyperref}
\usepackage{comment}
\usepackage[caption=false]{subfig}
\usepackage{multirow}

\providecommand{\U}[1]{\protect\rule{.1in}{.1in}}

\newtheorem{remark}[theorem]{Remark}
\begin{document}

\title{Inexact methods for symmetric stochastic eigenvalue problems\thanks{%
This work is based upon work supported by the U.\thinspace\ S.\thinspace\
Department of Energy Office of Advanced Scientific Computing Research,
Applied Mathematics program under Award Number DE-SC0009301, and by the
U.\thinspace\ S.\thinspace\ National Science Foundation under grant DMS1521563.}}
\author{Kookjin Lee\thanks{This work was performed while pursuing a Ph.D. in the Department of Computer Science at the University of Maryland, College Park, MD~20742. \textit{Current affiliation}:
Extreme-scale Data Science and Analytics Department, Sandia National Laboratories, Livermore, CA~94550 (\texttt{koolee@sandia.gov}).} \and Bed\v{r}ich Soused\'{\i}k\thanks{%
Department of Mathematics and Statistics, University of Maryland, Baltimore
County, 1000 Hilltop Circle, Baltimore, MD~21250 (\texttt{sousedik@umbc.edu}%
).} }

%Department of Computer Science, University of Maryland, College Park, MD 20742 (\texttt{klee@cs.umd.edu}) }

\maketitle

\begin{abstract}
We study two inexact methods for solutions of random eigenvalue problems in
the context of spectral stochastic finite elements. In particular, given a
parameter-dependent, symmetric matrix operator, the methods solve for
eigenvalues and eigenvectors represented using polynomial chaos expansions.
Both methods are based on the stochastic Galerkin formulation of the
eigenvalue problem and they exploit its Kronecker-product structure. 
The first method is an inexact variant of the stochastic inverse
subspace iteration [B. Soused\'{\i}k, H. C. Elman, \textit{SIAM/ASA Journal
on Uncertainty Quantification} 4(1), pp. 163--189, 2016]. The second method
is based on an inexact variant of Newton iteration. In both cases, the
problems are formulated so that the associated stochastic Galerkin matrices
are symmetric, and the corresponding linear problems are solved using
preconditioned Krylov subspace methods with several novel hierarchical
preconditioners. The accuracy of the methods is compared with that of Monte
Carlo and stochastic collocation, and the effectiveness of the methods is
illustrated by numerical experiments.
\end{abstract}



\pagestyle{myheadings}
\thispagestyle{plain}
\markboth{K. LEE AND B.\ SOUSED\'{I}K}{INEXACT METHODS FOR STOCHASTIC EIGENVALUE PROBLEMS}

\section{Introduction}

\label{sec:introduction}Eigenvalue analysis is important in a number of
applications, for example in modeling of vibrations of mechanical
structures, neutron transport criticality computations, or stability of
dynamical systems, to name a few. The behavior of the underlying
mathematical models depends on proper choice of parameters entering the
model through coefficients, boundary conditions or forces. However, in
practice the exact values of these parameters are not known and they are
treated as random processes. The uncertainty is translated by discretization
into the matrix operators and subsequently into eigenvalues and
eigenvectors. The standard techniques to solve this problems include Monte
Carlo~(MC) methods~\cite%
{Brockway-1985-MCE,Nightingale-2007-MCE,Pradlwarter-2002-REP}, which are
robust but relatively slow, and perturbation methods~\cite%
{Kaminski-2013-SPM,Kleiber-1992-SFE,Shinozuka-1972-REP,vomScheidt-1983-REP},
which are limited to models with low variability of uncertainty.

In this study, we use spectral stochastic finite element methods~(SSFEM)~%
\cite{Ghanem-1991-SFE,LeMaitre-2010-SMU,Lord-2014-ICS,Xiu-2010-NMS} for the
solution of symmetric eigenvalue problems. The assumption in these methods
is that the parametric uncertainty is described in terms of polynomials of
random variables, and they compute solutions that are also polynomials in
the same random variables in the so-called generalized polynomial chaos
(gPC) framework~\cite{Ghanem-1991-SFE,Xiu-2002-WAP}. There are two main
approaches: stochastic collocation~(SC) and stochastic Galerkin~(SG)
methods. The first approach is based on sampling, so the problem is
translated into a set of independent deterministic systems; the second one
is based on stochastic Galerkin projection and the problem is translated
into one large coupled deterministic system. While the SSFEM\ methods have
become quite popular for solving stochastic partial differential equations,
the literature addressing eigenvalue problems is relatively limited. The
stochastic inverse iteration in the context of the SG\ framework was
proposed by Verhoosel et al.~\cite{Verhoosel-2006-ISR}. Meidani and Ghanem~%
\cite{Meidani-2012-SMD,Meidani-2014-SPI} formulated stochastic subspace
iteration using a stochastic version of the modified Gram-Schmidt algorithm.
Soused\'{\i}k and Elman~\cite{Sousedik-2016-ISI} introduced stochastic
inverse subspace iteration by combining the two techniques, they showed that
deflation of the mean matrix can be used to compute expansions of the
interior eigenvalues, and they also showed that the stochastic Rayleigh
quotient alone provides a good approximation of an eigenvalue expansion; see
also~\cite{Chen-2017-SSE,Elman-2018-LRS,Pagnacco-2016-SIP} for closely related methods. The
authors of~\cite{Meidani-2014-SPI,Sousedik-2016-ISI} used a quadrature-based
normalization of eigenvectors. Normalization based on a solution of a small
nonlinear problem was proposed by Hakula et al.~\cite{Hakula-2015-AMS}, and
Hakula and Laaksonen~\cite{Hakula-2017-ACS} also provided an asymptotic
convergence theory for the stochastic iteration. In an alternative approach,
Ghanem and Ghosh~\cite{Ghanem-2007-ECR,Ghosh-2008-SCA} proposed two
numerical schemes---one based on Newton iteration and another based on an
optimization problem (see also~\cite{Ghosh-2013-ARE,Ghosh-2012-ISA}). Most
recently, Benner et al.~\cite{Benner-2017-INK} formulated an inexact
low-rank Newton--Krylov method, in which the stochastic Galerkin linear
systems are solved using the BiCGStab method with a variant of mean-based
preconditioner. In alternative approaches, Pascual and Adhikari~\cite%
{Pascual-2012-HPP} introduced several hybrid perturbation-polynomial chaos
methods, and Williams~\cite%
{Williams-2010-MSSnote,Williams-2010-MSS,Williams-2013-MSS} presented a
method that avoids the nonlinear terms in the conventional method of
stochastic eigenvalue calculation but introduces an additional independent
variable.

We formulate two inexact methods for symmetric eigenvalue problems
formulated in the SSFEM\ framework and based on the SG\ formulation. The
first method is an inexact variant of the stochastic inverse subspace
iteration from~\cite{Sousedik-2016-ISI}, in which the linear stochastic
Galerkin systems are solved using the conjugate gradient method with the
truncated hierarchical preconditioner~\cite{Sousedik-2014-THP} (see also~%
\cite{Sousedik-2014-HSC}). The second method is an inexact variant of the
Newton iteration from~\cite{Ghanem-2007-ECR}, in which the linear stochastic
Galerkin systems are solved using preconditioned MINRES\ and GMRES. The
methods are derived using the Kronecker-product formulation and we also
comment on the so-called matricized format. The formulation of the Newton's
method is closely related to that of~\cite{Benner-2017-INK}, however we
consider general parametrization of stochastic coefficients, the Jacobian
matrices are symmetrized, %using a simple algebraic trick, 
and we propose a class of hierarchical preconditioners, which can be viewed
as extensions of the hierarchical preconditioners used for the first method.
We also note that we have recently successfully combined an inexact
Newton--Krylov method with the stochastic Galerkin framework in a different
context~\cite{Lee-2017-LRS,Sousedik-2016-SGM}. The performance of both
methods is illustrated by numerical experiments, and the results are
compared to that of %Monte Carlo and stochastic collocation
MC\ and SC methods.

The paper is organized as follows. In Section~\ref{sec:problem} we introduce
the stochastic eigenvalue problem, in Section~\ref{sec:sampling} we recall
the solution techniques using sampling methods (Monte Carlo and stochastic
collocation), in Section~\ref{sec:sg} we introduce the stochastic Galerkin
formulation, in Section~\ref{sec:isisi}\ we formulate the inverse subspace
iteration and in Section~\ref{sec:ni} the Newton iteration,\ in
Section~\ref{sec:numerical}\ we report the results of numerical experiments,
and in Section~\ref{sec:conclusion}\ we summarize %and conclude 
our work.
In Appendices~\ref{sec:ini} and~\ref{sec:mvp} we describe algorithmic details, 
and in Appendix~\ref{sec:comp_cost} we discuss the computational cost.

\section{Stochastic eigenvalue problem}

\label{sec:problem}Let $D$ be a bounded physical domain, and let $\left(
\Omega ,\mathcal{F},\mathcal{P}\right) $ be a complete probability space,
that is, $\Omega $ is the {a} sample space with a $\sigma $-algebra$~%
\mathcal{F}$ and {a} probability measure $\mathcal{P}$. We assume that the
randomness in the mathematical model is induced by a vector $\xi :\Omega
\mapsto \Gamma \subset 
%TCIMACRO{\U{211d} }%
%BeginExpansion
\mathbb{R}
%EndExpansion
^{m_{\xi }}$\ of independent, identically distributed random variables $\xi
_{1}(\omega ),\dots ,\xi _{m_{\xi }}(\omega )$, where $\omega \in \Omega $.
Let $\mathcal{B}(\Gamma )$ denote the Borel $\sigma $-algebra on$~\Gamma $
induced by~$\xi $ and $\mu $ denote the induced measure. The expected value
of the product of measurable functions on$~\Gamma $ determines a Hilbert
space $T_{\Gamma }\equiv L^{2}\left( \Gamma ,\mathcal{B}(\Gamma ),\mu
\right) $ with inner product 
\begin{equation}
\left\langle u,v\right\rangle =\mathbb{E}\left[ uv\right] =\int_{\Gamma
}u(\xi )v(\xi )\,\mu (\xi )d\xi ,  \label{eq:stoch-inner-prod}
\end{equation}%
where the symbol $\mathbb{E}$ denotes the mathematical expectation.

In computations we will use a finite-dimensional subspace $T_{p}\subset
T_{\Gamma }$ spanned by a set of multivariate polynomials $\left\{ \psi
_{\ell }(\xi )\right\} $ that are orthonormal with respect to the density
function $\mu $, that is $\mathbb{E}\left[ \psi _{k}\psi _{\ell }\right]
=\delta _{k\ell }$, where$~\delta _{k\ell }$ is the Kronecker delta, and $%
\psi _{0}$ is constant. This will be referred to as the gPC basis~\cite%
{Xiu-2002-WAP}. The dimension of the space$~T_{p}$ depends on the polynomial
degree. For polynomials of total degree$~p$, the dimension is $n_{\xi }=%
\binom{m_{\xi }+p}{p}$. We suppose we are given a symmetric matrix-valued
random variable$~A(x,\xi )$ represented as%
\begin{equation}
A(x,\xi )=\sum_{\ell =1}^{n_{a}}A_{\ell }(x)\psi _{\ell }(\xi ),
\label{eq:stoch-exp-A}
\end{equation}%
where each$~A_{\ell }$ is a deterministic matrix of size $n_{x}\times n_{x}$
with size determined by the discretization of the physical domain, and $%
A_{1} $ is the mean value matrix, that is $A_{1}=\mathbb{E}\left[ A(x,\cdot )%
\right] $. The representation~(\ref{eq:stoch-exp-A}) is obtained from either
the Karhunen-Lo\`{e}ve expansion or, more generally, a stochastic expansion
of an underlying random process. 
%specific examples are provided in section with numerical experiments. 

We are interested in a solution of the following stochastic eigenvalue
problem: find a set of stochastic eigenvalues~$\lambda ^{s}$ and
corresponding eigenvectors $u^{s}$, $s=1,\dots ,n_{s}$, which almost surely
(a.s.) satisfy the equation 
\begin{equation}
A(x,\xi )u^{s}(x,\xi )=\lambda ^{s}(\xi )u^{s}(x,\xi ),\qquad \forall x\in D,
\label{eq:param_eig}
\end{equation}%
where $\lambda ^{s}(\xi )\in \mathbb{R}$ and $u^{s}(\xi )\in \mathbb{R}%
^{n_{x}}$, along with a normalization condition 
\begin{equation}
\left\langle {u^{s}(x,\xi )},u^{s}(x,\xi )\right\rangle _{%
%TCIMACRO{\U{211d} }%
%BeginExpansion
\mathbb{R}
%EndExpansion
}=1,  \label{eq:param_normal}
\end{equation}%
where $\left\langle \cdot ,\cdot \right\rangle _{%
%TCIMACRO{\U{211d} }%
%BeginExpansion
\mathbb{R}
%EndExpansion
}$ denotes the inner product of two vectors.

We will search for expansions of\ eigenpairs$~{(\lambda ^{s},u^{s})}$, $%
s=1,\dots ,n_{s}$, in the form 
\begin{equation}
\lambda ^{s}(\xi )=\sum_{k=1}^{n_{\xi }}\lambda _{k}^{s}\psi _{k}(\xi
),\qquad u^{s}(x,\xi )=\sum_{k=1}^{n_{\xi }}u_{k}^{s}\psi _{k}(\xi ),
\label{eq:sol_mat}
\end{equation}%
where $\lambda _{k}^{s}\in \mathbb{R}$ and $u_{k}^{s}\in \mathbb{R}^{n_{x}}$
are the coefficients corresponding to the basis $\left\{ \psi _{k}\right\} $%
. Equivalently to~(\ref{eq:sol_mat}), using the symbol~$\otimes $ for the
Kronecker product, we write 
\begin{equation}
\lambda ^{s}(\xi )=\Psi (\xi )^{T}\bar{\lambda}^{s},\qquad u^{s}(x,\xi
)=(\Psi (\xi )^{T}\otimes I_{n_{x}})\bar{u}^{s},  \label{eq:param_sol}
\end{equation}%
where $\Psi (\xi )=[\psi _{1}(\xi ),\ldots ,\psi _{n_{\xi }}(\xi )]^{T}$, $%
\bar{\lambda}^{s}=[\lambda _{1}^{s},\ldots ,\lambda _{n_{\xi }}^{s}]^{T}$,
and $\bar{u}^{s}=[(u_{1}^{s})^{T},\ldots (u_{n_{\xi }}^{s})^{T}]^{T}$. 
%We will explore and compare several ways to approximate these quantities. 
%and use Monte Carlo simulation to validate our approach.

\begin{remark}
One can in general consider different number of terms in the two expansions~(%
\ref{eq:sol_mat}). However, since the numerical experiments in~\cite%
{Sousedik-2016-ISI} and also in the present work indicate virtually no
effect when the number of terms in eigenvalue expansion is larger than in
the eigenvector expansion, we consider here the same number of terms in both
expansions, see also Remark~\ref{rem:RQ}.
\end{remark}

\subsection{Sampling methods}

\label{sec:sampling}Both Monte Carlo and stochastic collocation methods are
based on sampling. The coefficients are defined by a discrete projection 
\begin{equation}
\lambda _{k}^{s}=\left\langle \lambda ^{s},\psi _{k}\right\rangle ,\quad
k=1,\dots ,n_{\xi },\qquad u_{k}^{s}=\left\langle u^{s},\psi
_{k}\right\rangle ,\quad k=1,\dots ,n_{\xi }.  \label{eq:eig-gPC-proj}
\end{equation}%
The evaluations of coefficients in~(\ref{eq:eig-gPC-proj}) entail solving a
set of independent deterministic eigenvalue problems at a set of sample
points$~\xi ^{(q)}$, $q=1,\dots ,n_{MC}$ or$~n_{q}$, 
\begin{equation*}
A(\xi ^{(q)})u^{s}(\xi ^{(q)})=\lambda ^{s}\left( \xi ^{(q)}\right)
u^{s}\left( \xi ^{(q)}\right) ,\qquad s=1,\dots ,n_{s}.
\end{equation*}%
In the Monte Carlo method, the sample points$~\xi ^{(q)}$, $q=1,\dots
,n_{MC},$\ are generated randomly following the distribution of the random
variables$~\xi _{i}$, $i=1,\dots ,m_{\xi }$, and moments of solution are
computed by ensemble averaging. In addition, the coefficients in~(\ref%
{eq:sol_mat}) can be computed as\footnote{%
In numerical experiments we avoid projections on the gPC and work %directly
with the sampled quantities.}%, unless stated otherwise.} 
\begin{equation*}
\lambda _{k}^{s}=\frac{1}{n_{MC}}\sum_{q=1}^{n_{MC}}\lambda ^{s}(\xi
^{(q)})\psi _{k}\left( \xi ^{(q)}\right) ,\qquad u_{mk}^{s}=\frac{1}{n_{MC}}%
\sum_{q=1}^{n_{MC}}u^{s}(x_{m},\mathbb{\xi }^{(q)})\psi _{k}(\mathbb{\xi }%
^{(q)}),
\end{equation*}%
where $u_{mk}^{s}$ is the $m$th element of $u_{k}^{s}$. For stochastic
collocation%\footnote{%
%The stochastic collocation is used here in the form of so-called
%nonintrusive stochastic Galerkin method.}, 
, which is used here in the form of so-called nonintrusive stochastic
Galerkin method, the sample points$~\xi ^{(q)}$, $q=1,\dots ,n_{q},$ consist
of a predetermined set of \emph{collocation points}, and the coefficients $%
\lambda _{k}^{s}$ and $u_{k}^{s}$\ in expansions~(\ref{eq:sol_mat}) are
determined by evaluating~(\ref{eq:eig-gPC-proj}) in the sense of~(\ref%
{eq:stoch-inner-prod}) using numerical quadrature 
\begin{equation}
\lambda _{k}^{s}=\sum_{q=1}^{n_{q}}\lambda ^{s}(\mathbb{\xi }^{(q)})\psi
_{k}(\mathbb{\xi }^{(q)})w^{(q)},\qquad
u_{mk}^{s}=\sum_{q=1}^{n_{q}}u^{s}(x_{m},\mathbb{\xi }^{(q)})\psi _{k}(%
\mathbb{\xi }^{(q)})w^{(q)},  \label{eq:Q-lambda-u}
\end{equation}%
where $\mathbb{\xi }^{(q)}$ are the quadrature (collocation) points and $%
w^{(q)}$ are quadrature weights. We refer, e.g., to~\cite{LeMaitre-2010-SMU}%
\ for a discussion of quadrature rules. Details of the rule we use in our
numerical experiments are discussed in Section~\ref{sec:numerical}.

\subsection{Stochastic Galerkin formulation}

\label{sec:sg}The main contribution of this paper is the development of two
inexact methods based on the stochastic Galerkin formulation of eigenvalue
problem~(\ref{eq:param_eig})--(\ref{eq:param_normal}). The formulation
entails a projection 
\begin{align}
\left\langle Au^{s},\psi _{k}\right\rangle & =\left\langle \lambda
^{s}u^{s},\psi _{k}\right\rangle , & k& =1,\dots ,n_{\xi },\quad s=1,\dots
,n_{s},  \label{eq:SG-eig} \\
\left\langle u^{sT}u^{s},\psi _{k}\right\rangle & =\delta _{k1}, & k&
=1,\dots ,n_{\xi },\quad s=1,\dots ,n_{s}.  \label{eq:SG-normal}
\end{align}%
Let us introduce the notation %\begin{equation*}
%H_{\ell }=\left[ h_{\ell ,kj}\right] ,\quad h_{\ell ,kj}\equiv \mathbb{E}%
%\left[ \psi _{\ell }\psi _{k}\psi _{j}\right] ,\qquad \ell =1,\dots
%,n_{a},\quad j,k=1,\dots ,n_{\xi }.
%\end{equation*}%
\begin{equation}
\lbrack H_{\ell }]_{kj}=h_{\ell ,kj},\quad h_{\ell ,kj}\equiv \mathbb{E}%
\left[ \psi _{\ell }\psi _{k}\psi _{j}\right] ,\qquad \ell =1,\dots
,n_{a},\quad j,k=1,\dots ,n_{\xi }.  \label{eq:cijk}
\end{equation}%
Substituting~(\ref{eq:stoch-exp-A}) and~(\ref{eq:sol_mat}) into~(\ref%
{eq:SG-eig})--(\ref{eq:SG-normal}) yields a nonlinear system, 
%written in the tensor-product format as
\begin{eqnarray}
\left( \sum_{\ell =1}^{n_{a}}H_{\ell }\otimes A_{\ell }\right) \overline{u}%
^{s} &=&\left( \sum_{i=1}^{n_{\xi }}H_{i}\otimes \lambda
_{i}^{s}I_{n_{x}}\right) \overline{u}^{s},\quad s=1,\dots ,n_{s},
\label{eq:SG-eig-2} \\
\sum_{j=1}^{n_{\xi }}\sum_{i=1}^{n_{\xi }}\left[ H_{k}\circ \left\langle
u_{i}^{s},u_{j}^{s}\right\rangle _{%
%TCIMACRO{\U{211d} }%
%BeginExpansion
\mathbb{R}
%EndExpansion
}\right] _{ij} &=&\delta _{k1},\qquad k=1,\dots ,n_{\xi },\quad s=1,\dots
,n_{s},  \label{eq:SG-normal-2}
\end{eqnarray}%
where the symbol$~\circ $ is the Hadamard product, see, e.g.,~\cite[Chapter~5%
]{Horn-1991-TMA}. An equivalent formulation of~(\ref{eq:SG-eig-2})--(\ref%
{eq:SG-normal-2}) is obtained as follows. Substituting~(\ref{eq:param_sol})
into~(\ref{eq:param_eig})--(\ref{eq:param_normal}) and rearranging, we get 
\begin{eqnarray*}
(\Psi ^{T}\otimes A)\bar{u}^{s} &=&((\bar{\lambda}^{s})^{T}\Psi \Psi
^{T}\otimes I_{n_{x}})\bar{u}^{s}, \\
\bar{u}^{s}{}^{T}(\Psi (\xi )\Psi (\xi )^{T}\otimes I_{n_{x}})\bar{u}^{s}
&=&1,
\end{eqnarray*}%
and employing Galerkin projection~(\ref{eq:SG-eig})--(\ref{eq:SG-normal})
yields the equivalent formulation 
\begin{eqnarray}
\mathbb{E}[\Psi \Psi ^{T}\otimes A]\bar{u}^{s} &=&\mathbb{E}[((\bar{\lambda}%
^{s})^{T}\Psi )\Psi \Psi ^{T}\otimes I_{n_{x}})]\bar{u}^{s},
\label{eq:SG-eig-3} \\
\mathbb{E}[\Psi \otimes (\bar{u}^{s}{}^{T}(\Psi \Psi ^{T}\otimes I_{n_{x}})%
\bar{u}^{s})] &=&\mathbb{E}[\Psi \otimes 1].  \label{eq:SG-normal-3}
\end{eqnarray}

Finally, we note that the methods can be equivalently formulated in the
so-called \emph{matricized} format, which can also simplify the
implementation. To this end, we make use of isomorphism between~$%
%TCIMACRO{\U{211d} }%
%BeginExpansion
\mathbb{R}
%EndExpansion
^{n_{x}n_{\xi }}$ and $%
%TCIMACRO{\U{211d} }%
%BeginExpansion
\mathbb{R}
%EndExpansion
^{n_{x}\times n_{\xi }}$\ determined by the operators $\text{vec}$ and $%
\text{mat}$: $\bar{u}^{s}=\text{vec}(\bar{U}^{s})$, $\bar{U}^{s}=\text{mat}$(%
$\bar{u}^{s})$, where $\bar{u}^{s}\in 
%TCIMACRO{\U{211d} }%
%BeginExpansion
\mathbb{R}
%EndExpansion
^{n_{x}n_{\xi }}$, $\bar{U}^{s}\in 
%TCIMACRO{\U{211d} }%
%BeginExpansion
\mathbb{R}
%EndExpansion
^{n_{x}\times n_{\xi }}$ and the upper/lower case notation is assumed
throughout the paper, so $\bar{R}^{s}=\text{mat}$($\bar{r}^{s})$, etc.
Specifically, we define the \emph{matricized}\ coefficients of the
eigenvector expansion 
\begin{equation}
\bar{U}^{s}=\text{mat}(\bar{u}^{s})=\left[ u_{1}^{s},u_{2}^{s},\ldots
,u_{n_{\xi }}^{s}\right] \in 
%TCIMACRO{\U{211d} }%
%BeginExpansion
\mathbb{R}
%EndExpansion
^{n_{x}\times n_{\xi }},  \label{eq:U}
\end{equation}%
where the column~$k$ contains the coefficients associated with the basis
function$~\psi _{k}$.

In the rest of the paper we explore two methods for solving the eigenvalue
problem~(\ref{eq:SG-eig-2})--(\ref{eq:SG-normal-2}), resp.~(\ref{eq:SG-eig-3}%
)--(\ref{eq:SG-normal-3}): the first is based on inverse subspace iteration
(Section~\ref{sec:isisi}), and the second one is based on Newton iteration
(Section~\ref{sec:ni}).

\section{Inexact stochastic inverse subspace iteration}

\label{sec:isisi} 
We formulate an inexact variant of the inverse subspace iteration from~\cite%
{Sousedik-2016-ISI} for the solution of~(\ref{eq:SG-eig-2})--(\ref%
{eq:SG-normal-2}). 
%While inexact inverse iteration for deterministic problems was studied, 
%for example, in~\cite{Golub-2000-III,Xue-2009-CAI}, we
%note that we are not simply searching for eigenvalues of the stochastic
%Galerkin matrix given by the operator on the left-hand side of~(\ref%
%{eq:SG-eig-2}). 
Stochastic inverse iteration was formulated in~\cite{Verhoosel-2006-ISR} for
the case when a stochastic expansion of a single eigenvalue is sought. It
was suggested in~\cite{Sousedik-2016-ISI} that the matrix~$A_{1}$ can be
deflated, rather than applying a shift, to find an expansion of an interior
eigenvalue, and a stochastic version of modified Gram-Schmidt process~\cite%
{Meidani-2014-SPI} can be applied if more eigenvalues are of interest. 
%In this section, 
In this section, we formulate an inexact variant of the stochastic inverse
subspace iteration~\cite[Algorithm~3.2]{Sousedik-2016-ISI}, whereby the
linear systems~(\ref{eq:alg-SISI-solve}) are solved only approximately using
preconditioned conjugate gradient method~(PCG). %with the mean-based~\cite%
%{Pellissetti-2000-ISS,Powell-2009-BDP} and the truncated hierarchical
%Gauss-Seidel~\cite{Sousedik-2014-THP} preconditioners. 
The method is formulated as Algorithm~\ref{alg:isisi}. We now describe its
components in detail, % using the Kronecker-product notation.
and for simplicity we drop the superscript~$^{(n)}$ in the description.

\begin{algorithm}[hbpt]
\caption{Inexact stochastic inverse subspace iteration}
\label{alg:isisi}
\begin{algorithmic}[1]
\State Find the~$n_{s}$ smallest eigenpairs of  
\begin{align}
A_{1}\,w^s = \mu^s\, w^s, \qquad s=1,\dots,n_s. \label{eq:alg-SISI-ini}
\end{align}
\State \textbf{if} $\mu^1=\min (\mu^s)>0$, set $\rho=0$, 
            \textbf{else} shift $A_1=A_1+\rho I_{n_x}$, where $\rho>| \mu^1 |$. \textbf{end~if}
\State Initialize 
\begin{align}
u_{1}^{s,\left( 0\right) } = w^s,\quad u_{i}^{s,\left( 0\right) } = 0, \qquad 
s=1,\dots,n_{s},\quad i=2,\dots,n_{\xi}.
\end{align}
\For{$n=0,1,2,\dots$} 
\State Use conjugate gradients with preconditioner from Algorithm~\ref{alg:MB} or~\ref{alg:hGS} to solve %the stochastic Galerkin system 
\begin{equation}
\left( \sum_{\ell =1}^{n_{a}}H_{\ell }\otimes A_{\ell }\right) \overline{v}%
^{s,(n)} = \overline{u}^{s,(n)},\quad s=1,\dots ,n_{s}.  \label{eq:alg-SISI-solve}
\end{equation}
\If{$n_{s}=1$} normalize using the quadrature rule~(\ref{eq:vector-normalize}): 
$ \overline{u}^{1,\left( n+1\right) } \leftarrow\overline{v}^{1,\left( n\right) } $.
\Else{} orthogonalize using the stochastic modified Gram-Schmidt process: \newline
${} \qquad \qquad \overline{u}^{s,\left( n+1\right) } \leftarrow \overline{v}^{s,\left( n\right) }, \quad s=1,\dots,n_{s}$.
\EndIf
\State Check convergence.
\EndFor
\State Use the stochastic Rayleigh quotient~(\ref{eq:RQ}) to compute the eigenvalue expansions.
\State \textbf{if} $\rho>0$, shift $\lambda_1^s=\lambda_1^s-\rho$ for $s=1,\dots,n_s$. \textbf{end if}  
\end{algorithmic}
\end{algorithm}

\paragraph{Matrix-vector product}

The conjugate gradient method and computation of the stochastic Rayleigh
quotient require a stochastic version of a matrix-vector product, which
corresponds to evaluation of the projection 
\begin{equation*}
v_{k}^{s}=\left\langle v^{s},\psi _{k}\right\rangle =\left\langle
Au^{s},\psi _{k}\right\rangle ,\qquad k=1,\dots ,n_{\xi }.
\end{equation*}%
Since $(V\otimes W)\text{vec}(X)=\text{vec}(WXV^{T})$, the coefficients of
the expansion are 
\begin{equation}
\bar{v}^{s}=\mathbb{E}[\Psi \Psi ^{T}\otimes A]\bar{u}^{s}=\sum_{\ell
=1}^{n_{a}}({H}_{\ell }\otimes A_{\ell })\bar{u}^{s}\quad \Leftrightarrow
\quad \bar{V}^{s}=\sum_{\ell =1}^{n_{a}}A_{\ell }\bar{U}^{s}H_{\ell }^{T}.
\label{eq:mat-vec}
\end{equation}%
The use of this computation for the Rayleigh quotient is described below. We
also note that Algorithm~\ref{alg:isisi} can be modified to perform subspace
iteration~\cite[Algorithm~$4$]{Meidani-2014-SPI} for identifying the largest
eigenpairs. In this case, the solve~(\ref{eq:alg-SISI-solve}) is simply
replaced by a matrix-vector product~(\ref{eq:mat-vec}).

\paragraph{Stochastic Rayleigh quotient}

In the deterministic case, the Rayleigh quotient\ is used to compute the
eigenvalue corresponding to a normalized eigenvector$~u$ as $\lambda =u^{T}v$%
, where $v=Au$. For the stochastic Galerkin method, the Rayleigh quotient
defines the coefficients of a stochastic expansion of the eigenvalue defined
via a projection 
\begin{equation*}
\lambda _{k}^{s}=\left\langle \lambda ^{s},\psi _{k}\right\rangle
=\left\langle u^{s}{}^{T}v^{s},\psi _{k}\right\rangle ,\qquad k=1,\dots
,n_{\xi }.
\end{equation*}%
The coefficients of$~v^{s}$ are computed using~(\ref{eq:mat-vec}) and the
coefficients$~\lambda _{k}^{s}$ are 
\begin{equation*}
\lambda _{k}^{s}=\mathbb{E}\left[ \left( \left( \Psi ^{T}\otimes 1\right) 
\bar{\lambda}^{s}\right) \psi _{k}\right] =\mathbb{E}\left[ \left( \bar{u}%
^{sT}\left( \Psi \Psi ^{T}\otimes I_{n_{x}}\right) \bar{v}^{s}\right) \,\psi
_{k}\right] ,\qquad k=1,\dots ,n_{\xi },
\end{equation*}%
which is 
\begin{equation}
\lambda _{k}^{s}=\sum_{j=1}^{n_{\xi }}\sum_{i=1}^{n_{\xi }}\left[ H_{k}\circ
\left\langle u_{i}^{s},v_{j}^{s}\right\rangle _{%
%TCIMACRO{\U{211d} }%
%BeginExpansion
\mathbb{R}
%EndExpansion
}\right] _{ij}=\sum_{j=1}^{n_{\xi }}\sum_{i=1}^{n_{\xi }}\left[ H_{k}\circ
\left( \bar{U}^{sT}\bar{V}^{s}\right) \right] _{ij},\qquad k=1,\dots ,n_{\xi
}.  \label{eq:RQ}
\end{equation}

\begin{remark}
\label{rem:RQ}The Rayleigh quotient~(\ref{eq:RQ}) finds~$n_{\xi }$
coefficients of the eigenvalue expansion, which is consistent with Newton
iteration formulated in Section~\ref{sec:ni} and also with the literature~%
\cite{Meidani-2014-SPI,Verhoosel-2006-ISR}. We note that it would be
possible to compute the coefficients$~\lambda _{k}$ for $k>n_{\xi }$ as
well, because the inner product~$u^{T}v$ of two eigenvectors which are
expanded using chaos polynomials up to degree~$p$ has nonzero chaos
coefficients up to degree~$2p$. An alternative is to use a \emph{full}
representation of the Rayleigh quotient based on the projection of $u^{T}Au$%
. However, from our experience in the present and the previous work~\cite%
{Sousedik-2016-ISI}, the representation~(\ref{eq:RQ}) is sufficient.
\end{remark}

\paragraph{Normalization and the Gram-Schmidt process}

Let $~\left\Vert \cdot \right\Vert _{2}$ denote the vector norm, induced by
the inner product$~\left\langle \cdot ,\cdot \right\rangle _{%
%TCIMACRO{\U{211d} }%
%BeginExpansion
\mathbb{R}
%EndExpansion
}$. That is, for a vector$~u$ evaluated at a point$~\xi $, 
\begin{equation}
\left\Vert u\left( \xi \right) \right\Vert _{2}=\sqrt{\sum_{n=1}^{n_{x}}%
\left( \left[ u\left( \xi \right) \right] _{n}\right) ^{2}}.
\label{eq:vector-norm}
\end{equation}
At each step of stochastic iteration the coefficients of a given set of
vectors$~\left\{ v^{s}\right\} _{s=1}^{n_{s}}$ are transformed into an
orthonormal set$~\left\{ u^{s}\right\} _{s=1}^{n_{s}}$ such that the
condition 
\begin{equation}
\left\langle u^{s}\left( \xi \right) ,u^{t}\left( \xi \right) \right\rangle
_{\mathbb{R}}=\delta _{st},\qquad \text{a.s}.,  \label{eq:orthonormal}
\end{equation}
and in particular~(\ref{eq:SG-normal-2}), is satisfied. We adopt the same
strategy as in$~$\cite{Meidani-2014-SPI,Sousedik-2016-ISI}, whereby the
coefficients of the orthonormal eigenvectors are calculated using a discrete
projection and a quadrature rule. An alternative approach to normalization,
based on solution of a relatively small nonlinear system was proposed by
Hakula et al.~\cite{Hakula-2015-AMS}.

Let us first consider \emph{normalization} of a vector, so $n_s=1$. The
coefficients in column$~k$ of$~\bar{U}^{1}$ corresponding to coefficients of
a normalized vector are computed as%
\begin{equation}
u_{k}^{1}=\sum_{q=1}^{n_{q}}\frac{v^{1}\left( \xi ^{\left( q\right) }\right) 
}{\left\Vert v^{1}\left( \xi ^{\left( q\right) }\right) \right\Vert _{2}}%
\,\psi _{k}\left( \xi ^{\left( q\right) }\right) \,w^{\left( q\right) }.
\label{eq:vector-normalize}
\end{equation}%
When $n_s>1$, the \emph{orthonormalization}~(\ref{eq:orthonormal}) is
performed by a combination of stochastic Galerkin projection and the
modified Gram-Schmidt algorithm as proposed in~\cite{Meidani-2014-SPI}, 
\begin{equation}
\mathbb{E}\left[ \Psi \otimes u^{s}\right] =\mathbb{E}\left[ \Psi \otimes
v^{s}\right] -\sum_{t=1}^{s-1}\mathbb{E}\left[ \Psi \otimes \left( \frac{%
\left\langle v^{s},u^{t}\right\rangle _{\mathbb{R}}}{\left\langle
u^{t},u^{t}\right\rangle _{\mathbb{R}}}u^{t}\right) \right] ,\quad s=2,\dots
,n_{s},  \label{eq:vector-orthogonalize}
\end{equation}%
Using the expansion~(\ref{eq:param_sol}) and rearranging, the coefficients
in column$~k$ of$~\bar{U}^{s}$\ are 
\begin{equation*}
u_{k}^{s}=v_{k}^{s}-\sum_{t=1}^{s-1}\chi _{k}^{ts},\qquad k=1,\dots ,n_{\xi
},\quad s=2,\dots ,n_{s},
\end{equation*}%
where 
\begin{equation*}
\chi ^{ts}(\xi )=\frac{\left\langle v^{s}(\xi ),u^{t}(\xi )\right\rangle _{%
%TCIMACRO{\U{211d} }%
%BeginExpansion
\mathbb{R}
%EndExpansion
}}{\left\langle u^{t}(\xi ),u^{t}(\xi )\right\rangle _{%
%TCIMACRO{\U{211d} }%
%BeginExpansion
\mathbb{R}
%EndExpansion
}}u^{t}(\xi ),
\end{equation*}
and the coefficients $\chi _{k}^{ts}$ are computed using a discrete
projection %and a quadrature rule
as in (\ref{eq:Q-lambda-u}), 
\begin{equation*}
\chi _{k}^{ts}=\sum_{q=1}^{n_{q}}\chi ^{ts}\left( \xi ^{\left( q\right)
}\right) \,\psi _{k}\left( \xi ^{\left( q\right) }\right) \,w^{\left(
q\right) }.
\end{equation*}

\paragraph{Stopping criteria}

The inexact iteration entails in each step of Algorithm~\ref{alg:isisi}\ a
solution of the stochastic Galerkin problem~(\ref{eq:alg-SISI-solve}) using
the preconditioned conjugate gradient method. We use the criteria proposed
by Golub and Ye~\cite[Eq.~(1)]{Golub-2000-III}; the criteria is satisfied
when the relative residual of PCG\ gets smaller than a factor of the
nonlinear residual from the previous step, that is 
\begin{equation}
\frac{\Vert \overline{u}^{s,(n)}\!-\!\left( \sum_{\ell =1}^{n_{a}}H_{\ell
}\otimes A_{\ell }\right) \!\overline{v}^{s,(n)}\Vert _{2}}{\Vert \overline{u%
}^{s,(n)}\Vert _{2}}\!<\!\tau \left\Vert {\left( \sum_{\ell
=1}^{n_{a}}H_{\ell }\otimes A_{\ell }\!-\!\sum_{i=1}^{n_{\xi }}H_{i}\otimes
\lambda _{i}^{s,(n-1)}I_{n_{x}}\right) \!\overline{u}^{s,(n-1)}}\right\Vert
_{2},  \label{eq:pcg-stop}
\end{equation}%
where the factor $\tau =10^{-2}$. 
%A similar condition is used in inexact Newton iteration introduced in the next section. 
It is important to note that Algorithm~\ref{alg:isisi} provides only the
coefficients of expansion of the projection of residual on the gPC\ basis,
that is 
\begin{equation}
\widetilde{r}_{k}^{s}=\left\langle Au^{s}-\lambda ^{s}u^{s},\psi
_{k}\right\rangle ,\qquad k=1,\dots ,n_{\xi },\quad s=1,\dots ,n_{s}
\label{eq:res-indicator}
\end{equation}%
One could assess accuracy using Monte Carlo sampling of this residual by
computing 
\begin{equation*}
r^{s}\left( \xi ^{i}\right) =A\left( \xi ^{i}\right) u^{s}\left( \xi
^{i}\right) -\lambda ^{s}\left( \xi ^{i}\right) u^{s}\left( \xi ^{i}\right)
,\qquad i=1,\dots ,N_{MC},\quad s=1,\dots ,n_{s}.
\end{equation*}%
However, in the numerical experiments we use a much less expensive
computation, which is based on using coefficients$~\widetilde{r}_{k}^{s}$
directly as an error indicator. In particular, we monitor the norms of the
terms of~$\widetilde{r}_{k}^{s}$ corresponding to expected value and
variance, 
\begin{equation}
\varepsilon _{1}^{s,\left( it\right) }=\left\Vert \widetilde{r}%
_{1}^{s,\left( n\right) }\right\Vert _{2},\qquad \varepsilon _{\sigma
^{2}}^{s,\left( it\right) }=\left\Vert \sum_{k=2}^{n_{\xi }}\left( 
\widetilde{r}_{k}^{s,\left( n\right) }\right) ^{2}\right\Vert _{2},\qquad
s=1,\dots ,n_{s}.  \label{eq:eps}
\end{equation}

\subsection{Preconditioners for the stochastic inverse iteration}

We use two preconditioners for problem~(\ref{eq:alg-SISI-solve}) -- the
mean-based preconditioner~\cite{Pellissetti-2000-ISS,Powell-2009-BDP} and
the hierarchical Gauss-Seidel preconditioner~\cite{Sousedik-2014-THP}. Both
preconditioners are formulated in the Kronecker-product format and we also
comment on the matricized formulation. We assume that a preconditioner$%
~M_{1} $ for the mean matrix$~A_{1}$ is available.

The mean-based preconditioner~(MB) is listed as Algorithm~\ref{alg:MB}.
Since $H_{1}=I_{n_{\xi }}$, the preconditioner entails$~n_{\xi }$\ block
diagonal solves with$~M_{1}$, and recalling that we can write $\bar{R}^{s}=%
\text{mat}(\bar{r}^{s})$, $\bar{V}^{s}=\text{mat}(\bar{v}^{s})$, its action
can be equivalently obtained by solving 
\begin{equation}
M_{1}\bar{V}^{s}=\bar{R}^{s}.  \label{eq:MB-matricized}
\end{equation}

\begin{algorithm}[hptb]
\caption{{\cite{Pellissetti-2000-ISS,Powell-2009-BDP}} Mean-based preconditioner (MB)}
\label{alg:MB} 
The preconditioner $M_{\text{MB}}: \bar{r}^{s} \longmapsto \bar{v}^{s}$
 for~(\ref{eq:alg-SISI-solve}) is defined as
\begin{equation*}
\left( H_1 \otimes M_1 \right)  \bar{v}^{s} = \bar{r}^{s}.  \label{eq:algMB1}
\end{equation*}
\end{algorithm}

\begin{algorithm}[hptb]
\caption{{\cite[Algorithm~3]{Sousedik-2014-THP}} Hierarchical Gauss-Seidel preconditioner (hGS)}
\label{alg:hGS} 
The preconditioner $M_{\text{hGS}}:\bar{r}^{s}\longmapsto \bar{v}^{s}$\ for~(%
\ref{eq:alg-SISI-solve}) is defined as follows.
\begin{algorithmic}[1]
\State Set the initial solution$~\bar{v}^{s}$ to zero and update in the following steps: 
\State Solve 
\begin{equation}
M_{1}v_{1}^{s}=r_{1}^{s}-\mathcal{F}_{1}v_{\left( 2:n_{\xi }\right) }^{s},
\qquad \text{ where }
\mathcal{F}_{1}=\sum_{t\in \mathcal{I}_{t}}\left( \left[ h_{t,\left(1\right) \left(
2:n_{\xi }\right) }\right] \otimes A_{t}\right). 
\label{eq:alg-hGS1}
\end{equation}

\For{$d=1,\ldots, p-1$} 
\State Set
$\ell =\left( n_{\ell }+1:n_{u}\right) ,\text{ where }n_{\ell }=\binom{m_{\xi
}+d-1}{d-1}\text{ and }n_{u}=\binom{m_{\xi }+d}{d}$.
\State Solve
\begin{equation}
\left( I_{n_{u}-n_{\ell }}\otimes M_{1}\right) v_{\left( \ell \right)
}^{s}=r_{\left( \ell \right) }^{s}-\mathcal{E}_{d+1}v_{\left( 1:n_{\ell
}\right) }^{s}-\mathcal{F}_{d+1}v_{\left( n_{u}+1:n_{\xi }\right) }^{s},
\label{eq:alg-hGS2}
\end{equation}%
where 
\begin{equation*}
\mathcal{E}_{d+1}=\sum_{t\in \mathcal{I}_{t}}\left( \left[ h_{t,\left( \ell
\right) \left( 1:n_{\ell }\right) }\right] \otimes A_{t}\right) ,\qquad 
\mathcal{F}_{d+1}=\sum_{t\in \mathcal{I}_{t}}\left( \left[ h_{t,\left( \ell
\right) \left( n_{u}+1:n_{\xi }\right) }\right] \otimes A_{t}\right). 
\end{equation*}
\EndFor

\State Set $\ell =\left( n_{u}+1:n_{\xi }\right)$.
\State Solve
\begin{equation*}
\left( I_{n_{\xi }-n_{u}}\otimes M_{1}\right) v_{\left( \ell \right)
}^{s}=r_{\left( \ell \right) }^{s}-\mathcal{E}_{p+1}v_{\left( 1:n_{u}\right)
}^{s},
\quad \text{ where } 
\mathcal{E}_{p+1}=\sum_{t\in \mathcal{I}_{t}}\left( \left[ h_{t,\left( \ell
\right) \left( 1:n_{u}\right) }\right] \otimes A_{t}\right). 
\end{equation*}

\For{$d=p-1,\ldots, 1$}
\State Set $\ell =\left( n_{\ell }+1:n_{u}\right) ,\text{ where }n_{\ell }=\binom{m_{\xi
}+d-1}{d-1}\text{ and }n_{u}=\binom{m_{\xi }+d}{d}$.
\State Solve~(\ref{eq:alg-hGS2}).
\EndFor
\State Solve~(\ref{eq:alg-hGS1}).
\end{algorithmic}
\end{algorithm}

The hierarchical Gauss-Seidel preconditioner~(hGS) is listed as Algorithm~%
\ref{alg:hGS}. We will denote by $~v_{\left( i:n\right) }^{s}$ a subvector
of~$\bar{v}^{s}$\ containing gPC coefficients~$i,i+1,\dots ,n$, and, in
particular,$\ \bar{v}^{s}=v_{\left( 1:n_{\xi }\right) }^{s}$. There are two
components of the preconditioner. The first component consists of
block-diagonal solves with blocks of varying sizes, but computed just as in
Algorithm~\ref{alg:MB}, resp. in~(\ref{eq:MB-matricized}). The second
component is used in the setup of the right-hand sides for the solves and
consists of matrix-vector products by certain subblocks of the stochastic
Galerkin matrix by vectors of corresponding sizes. To this end, we will
write $\left[ h_{t,(\ell )(k)}\right] $, with $(\ell )$\ and $(k)$ denoting
a set of (consecutive) rows and columns of matrix~$H_{t}$ so that, in
particular, ${H}_{t}=\left[ h_{t,(1:n_{\xi })(1:n_{\xi })}\right] $. Then,
the matrix-vector products can be written, cf.~(\ref{eq:mat-vec}) and note
the symmetry of$~H_{t}$, as 
\begin{equation}
v_{(\ell )}^{s}=\sum_{t\in \mathcal{I}_{t}}(\left[ h_{t,(\ell )(k)}\right]
\otimes A_{t})u_{(k)}^{s}\quad \Leftrightarrow \quad V_{(\ell
)}^{s}=\sum_{t\in \mathcal{I}_{t}}A_{t}U_{(k)}^{s}\left[ h_{t,(k)(\ell )}%
\right] ,  \label{eq:matvec-block}
\end{equation}%
where$~\mathcal{I}_{t}$\ is an index set$~\mathcal{I}_{t}\subseteq \left\{
1,\dots ,n_{\xi }\right\} $ indicating that the matrix-vector products may
be truncated. Possible strategies for truncation are discussed in~\cite%
{Sousedik-2014-THP}. In this study, we use $\mathcal{I}_{t}=\left\{ 1,\dots
,n_{t}\right\} $ with $n_{t}=\binom{m_{\xi }+p_{t}}{p_{t}}$ for some $%
p_{t}\leq p$ and, in particular, we set $t=\{0,1,2\}$. We also note that,
since the initial guess is zero in Algorithm~\ref{alg:hGS}, the
multiplications by$~\mathcal{F}_{1}$\ and $\mathcal{F}_{d+1}$ vanish from~(%
\ref{eq:alg-hGS1})--(\ref{eq:alg-hGS2}).

\section{Newton iteration}

\label{sec:ni}
\begin{figure}[!t]
\centering
\begin{tabular}{cc}
\includegraphics[width=6cm]{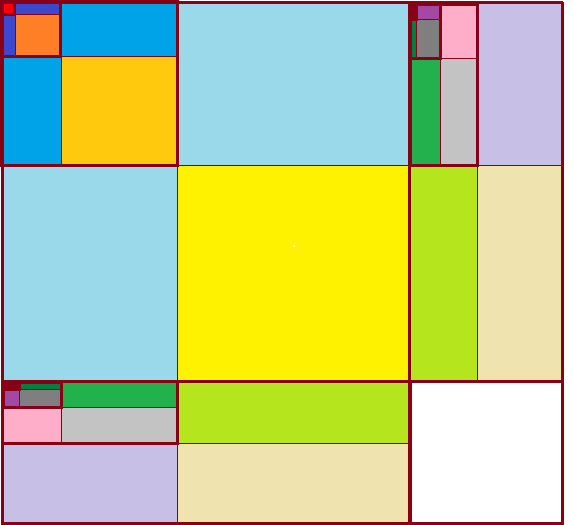} \hspace{3mm} %
\includegraphics[width=5.9cm]{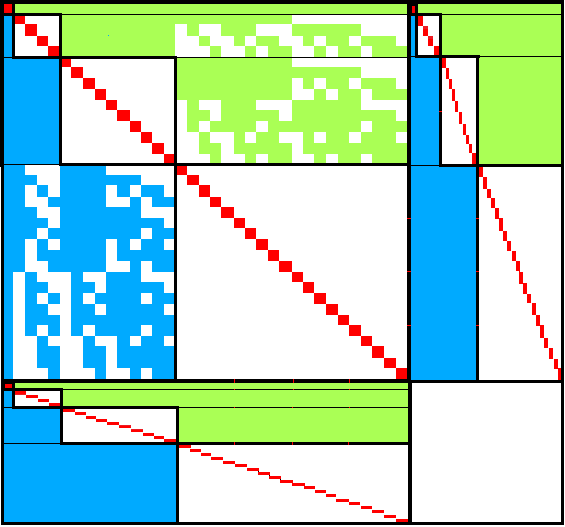} & 
\end{tabular}%
\caption{Hierarchical structure of the symmetric Jacobian matrix from~(%
\protect\ref{eq:Newton-mod}) (left) and splitting operator for the
constraint hierarchical Gauss-Seidel preconditioner from Algorithm~\protect
\ref{alg:chGS}--\protect\ref{alg:chGS_cont} (right).}
\label{fig:J-hierarchy}
\end{figure}

Use of Newton iteration to solve~(\ref{eq:SG-eig})--(\ref%
{eq:SG-normal}) was proposed in~\cite{Ghanem-2007-ECR}, and most recently
studied\ in~\cite{Benner-2017-INK}. We use a similar strategy also here and
formulate a line-search Newton method as Algorithm~\ref{alg:line_search}. To
begin, we consider the system of nonlinear equations~(\ref{eq:SG-eig-3})--(%
\ref{eq:SG-normal-3}) and rewrite it as 
\begin{equation}
\begin{bmatrix}
F(\bar{u}^{s},\bar{\lambda}^{s}) \\ 
G(\bar{u}^{s})%
\end{bmatrix}%
=0,\qquad s=1,\dots ,n_{s},  \label{eq:N-system}
\end{equation}%
where 
\begin{align}
F(\bar{u}^{s},\bar{\lambda}^{s})& \equiv \mathbb{E}[\Psi \Psi ^{T}\otimes A]%
\bar{u}^{s}-\mathbb{E}[((\bar{\lambda}^{s})^{T}\Psi )\Psi \Psi ^{T}\otimes
I_{n_{x}}]\bar{u}^{s},  \label{eq:F_mat} \\
G(\bar{u}^{s})& \equiv \mathbb{E}[\Psi \otimes \left( (\bar{u}%
^{s}{}^{T}(\Psi \Psi ^{T}\otimes I_{n_{x}})\bar{u}^{s})-1\right) ].
\label{eq:G_mat}
\end{align}%
The Jacobian matrix of$~$(\ref{eq:N-system}) is 
\begin{equation}
\mathcal{J}(\bar{u}^{s},\bar{\lambda}^{s})=%
\begin{bmatrix}
\frac{\partial F}{\partial \bar{u}^{s}} & \frac{\partial F}{\partial \bar{%
\lambda}^{s}} \\ 
\frac{\partial G}{\partial \bar{u}^{s}} & 0%
\end{bmatrix}%
,  \label{eq:jac}
\end{equation}%
where 
\begin{align}
\frac{\partial F}{\partial \bar{u}^{s}}(\bar{\lambda}^{s})& =\mathbb{E}[\Psi
\Psi ^{T}\otimes A]-\mathbb{E}[((\bar{\lambda}^{s})^{T}\Psi )\Psi \Psi
^{T}\otimes I_{n_{x}}],  \label{eq:jac_Fu} \\
\frac{\partial F}{\partial \bar{\lambda}^{s}}(\bar{u}^{s})& =-\mathbb{E}%
[\Psi ^{T}\otimes (\Psi \Psi ^{T}\otimes I_{n_{x}})\bar{u}^{s}],
\label{eq:jac_Fl} \\
\frac{\partial G}{\partial \bar{u}^{s}}(\bar{u}^{s})& =2\mathbb{E}[\Psi
\otimes ((\bar{u}^{s})^{T}(\Psi \Psi ^{T}\otimes I_{n_{x}}))].
\label{eq:jac_Gu}
\end{align}%
Step$~n$ of Newton iteration entails solving a linear system 
\begin{equation}
\left[ 
\begin{array}{cc}
\frac{\partial F}{\partial \bar{u}^{s}}(\bar{\lambda}^{s,(n)}) & \frac{%
\partial F}{\partial \bar{\lambda}^{s}}(\bar{u}^{s,(n)}) \\ 
\frac{\partial G}{\partial \bar{u}^{s}}(\bar{u}^{s,(n)}) & 0%
\end{array}%
\right] \left[ 
\begin{array}{c}
\delta \overline{u}^{s} \\ 
\delta \overline{\lambda }^{s}%
\end{array}%
\right] =-\left[ 
\begin{array}{c}
F(\bar{u}^{s,(n)},\bar{\lambda}^{s,(n)}) \\ 
G(\bar{u}^{s,(n)})%
\end{array}%
\right] ,  \label{eq:Newton}
\end{equation}%
followed by an update of the solution 
\begin{equation}
\left[ 
\begin{array}{c}
\overline{u}^{s,(n+1)} \\ 
\overline{\lambda }^{s,(n+1)}%
\end{array}%
\right] =\left[ 
\begin{array}{c}
\overline{u}^{s,(n)} \\ 
\overline{\lambda }^{s,(n)}%
\end{array}%
\right] +\left[ 
\begin{array}{c}
\delta \overline{u}^{s} \\ 
\delta \overline{\lambda }^{s}%
\end{array}%
\right] .  \label{eq:Newton-update}
\end{equation}%
The matrix$~\mathcal{J}(\bar{u}^{s},\bar{\lambda}^{s})$ is non-symmetric,
but since $\frac{\partial F}{\partial \bar{\lambda}^{s}}(\bar{u}^{s,(n)})=%
\left[ -\frac{1}{2}\frac{\partial G}{\partial \bar{u}^{s}}(\bar{u}^{s,(n)})%
\right] ^{T}$, we modify linear system~(\ref{eq:Newton}) in our
implementation as 
\begin{equation}
\left[ 
\begin{array}{cc}
\frac{\partial F}{\partial \overline{u}^{s}}(\bar{\lambda}^{s,(n)}) & \frac{%
\partial F}{\partial \bar{\lambda}^{s}}(\bar{u}^{s,(n)}) \\ 
\left[ \frac{\partial F}{\partial \bar{\lambda}^{s}}(\bar{u}^{s,(n)})\right]
^{T} & 0%
\end{array}%
\right] \left[ 
\begin{array}{c}
\delta \overline{u}^{s} \\ 
\delta \overline{\lambda }^{s}%
\end{array}%
\right] =\left[ 
\begin{array}{c}
-F(\bar{u}^{s,(n)},\bar{\lambda}^{s,(n)}) \\ 
\frac{1}{2}G(\bar{u}^{s,(n)})%
\end{array}%
\right] ,  \label{eq:Newton-mod}
\end{equation}%
which restores symmetry of linear systems solved in each step of Newton
iteration. The symmetric Jacobian matrix in~(\ref{eq:Newton-mod}) will be
denoted by$~J(\bar{u}^{s,(n)},\bar{\lambda}^{s,(n)})$. The hierarchical
structure of the Jacobian matrix, which is due to the stochastic Galerkin
projection, is illustrated by the left panel of Figure~\ref{fig:J-hierarchy}.
The systems~(\ref{eq:Newton-mod}) are solved inexactly using a preconditioned Krylov subspace method, 
and the details of evaluation of the right-hand side and the matrix-vector product 
are given in Appendix~\ref{sec:ini}.

\subsection{Inexact line-search Newton method}

\label{sec:line_search}In order to improve global convergence behavior of
Newton iteration, we consider a line-search modification of the method
following~\cite[Algorithm~11.4]{Nocedal-1999-NO}. To begin, let us define
the merit function as the sum of squares, 
\begin{equation*}
f(\bar{u}^{s,(n)},\bar{\lambda}^{s,(n)})=\frac{1}{2}\Vert r(\bar{u}^{s,(n)},%
\bar{\lambda}^{s,(n)})\Vert _{2}^{2},
\end{equation*}%
where $r$ is the residual of~(\ref{eq:N-system}), and denote 
\begin{equation*}
f_{n}=f(\bar{u}^{s,(n)},\bar{\lambda}^{s,(n)}),\qquad r_{n}=r(\bar{u}%
^{s,(n)},\bar{\lambda}^{s,(n)}),\qquad J_{n}=J(\bar{u}^{s,(n)},\bar{\lambda}%
^{s,(n)}).
\end{equation*}%
As the initial approximation of the solution, we use the eigenvectors and
eigenvalues of the associated mean problem given by the matrix $A_{1}$
concatenated by zeros, that is $\bar{u}^{s,(0)}=[(u_{1}^{s,(0)})^{T},0,%
\ldots ]^{T}$ and $\bar{\lambda}^{s,(0)}=[\lambda _{1}^{s,(0)},0,\ldots
]^{T} $, and the initial residual is 
\begin{equation*}
r_{0}=%
\begin{bmatrix}
F(\bar{u}^{s,(0)},\bar{\lambda}^{s,(0)}) \\ 
G(\bar{u}^{s,(0)})%
\end{bmatrix}%
.
\end{equation*}%
The line-search Newton method is summarized in our setting as Algorithm~\ref%
{alg:line_search}, and the choice of parameters$~\rho $ and$~c$ in the
numerical experiments is discussed in Section~\ref{sec:numerical}.

\begin{algorithm}[hptb]
\caption{\cite[Algorithm 11.4]{Nocedal-1999-NO} line-search Newton method}
\begin{algorithmic}[1]
\State Given $\rho,c \in (0, 1)$, set $\alpha^\ast = 1$.
\State Set $\bar u^{(0)}$ and $\bar \lambda^{(0)}$.
\For{$n=0,1,2,\ldots$}
\State $J_n p_n = - r_n$ \hfill (Find the Newton update $p_n$.) \label{ln:lin_sys}
\State $\begin{bmatrix} \delta \bar u^{(n)} \\ \delta \bar \lambda^{(n)} \end{bmatrix} = p_n$ 
\State $\alpha_n = \alpha^\ast$
\While{ $f(\bar u^{(n)} + \alpha_n \delta \bar u^{(n)}, \bar \lambda^{(n)} + \alpha_n \delta \bar \lambda^{(n)}) > f_n + c \, \alpha_n \nabla f_n^T p_n  $}
\State $\alpha_n \leftarrow \rho \, \alpha_n$
\EndWhile
\State $\bar u^{(n+1)} \leftarrow \bar u^{(n)} + \alpha_n \delta \bar u^{(n)}$ 
\State $\bar \lambda^{(n+1)} \leftarrow \bar \lambda^{(n)} + \alpha_n \delta \bar \lambda^{(n)}$
\State Check for convergence. \label{ln:ls_conv-check}
\EndFor
\end{algorithmic}
\label{alg:line_search}
\end{algorithm}

The inexact iteration entails in each step a solution of the stochastic
Galerkin linear system in Line~\ref{ln:lin_sys} of Algorithm~\ref%
{alg:line_search} given by~(\ref{eq:Newton-mod}) using a Krylov subspace
method. In our algorithm we use the adaptive stopping criteria for the
method, %(stopping tolerance \TEXTsymbol{<} tau * current\_residual,) 
\begin{equation}
\frac{\Vert r_{n}+J_{n}p_{n}\Vert _{2}}{\Vert r_{n}\Vert _{2}}<\tau
\left\Vert r{_{n-1}}\right\Vert _{2},  \label{eq:gmres-stop}
\end{equation}%
where $\tau =10^{-1}$. The for-loop is terminated when the convergence check
in Line~\ref{ln:ls_conv-check} is satisfied; in our numerical experiments we
check if $\left\Vert r_{n}\right\Vert _{2}<10^{-10}$.

\subsection{Preconditioners for the Newton iteration}

\label{sec:prec}The Jacobian matrices in~(\ref{eq:Newton-mod}) are
symmetric, indefinite, and so the linear systems can be ideally solved using 
\textsc{MINRES} iterative method. It is well known that a preconditioner for 
\textsc{MINRES} must be symmetric and positive definite cf., e.g.,~\cite%
{Wathen-2015-P}. A popular choice is a block diagonal preconditioner, cf.~%
\cite{Murphy-2000-NPI}, 
\begin{equation*}
\left[ 
\begin{array}{cc}
\widetilde{A} & 0 \\ 
0 & \widetilde{S}%
\end{array}%
\right] ,
\end{equation*}%
where~$\widetilde{A}\approx A$ and the Schur complement $\widetilde{S}%
\approx BA^{-1}B^{T}$ are obtained as approximations of the blocks in~(\ref%
{eq:jac-scheme}). 
%,cf.~(\ref{eq:jac_Fu}) and~(\ref{eq:jac_Gu}) scaled by $-1/2$. 
Such preconditioner, based on truncation of the series in$~$(\ref%
{eq:jac_Fu-impl}) and~(\ref{eq:jac_Fl-impl}) to the very first term, was
used in~\cite{Benner-2017-INK}. In such setup, we get 
\begin{eqnarray*}
\widetilde{A} &=&I_{n_{\xi }}\otimes A_{1}-(\lambda _{1}^{s}I_{n_{\xi
}}\otimes I_{n_{x}})=I_{n_{\xi }}\otimes (A_{1}-\lambda _{1}^{s}I_{n_{x}}) \\
&\approx &I_{n_{\xi }}\otimes (1-\lambda _{1}^{s})(A_{1}-I_{n_{x}}) \\
&\approx &I_{n_{\xi }}\otimes M_{1}^{s},
\end{eqnarray*}%
where the second line was used in~\cite{Benner-2017-INK}. In this study, we
use the third line with 
\begin{equation}
M_{1}^{s}=A_{1}-\epsilon _{M}\,\mu ^{s}I_{n_{x}},  \label{eq:M_1^s}
\end{equation}%
where $\mu ^{s}$ is the eigenvalue of the mean problem, cf.~(\ref%
{eq:alg-SISI-ini}). We note that it might be desirable to set the parameter $%
\epsilon _{M}\approx 1$, but $\epsilon _{M}\neq 1$ in order to guarantee
nonsingular$~M_{1}^{s}$, however more details for setup and use of~(\ref%
{eq:M_1^s}) are given in numerical experiments. Considering the first column
of~(\ref{eq:jac_Fl-impl}), cf.~(\ref{eq:jac_Fl}) and~(\ref{eq:jac_Fl-impl-2}%
), we get 
\begin{equation*}
\widetilde{B}^{T}=-(I_{n_{\xi }}\otimes u_{1}^{s}),
\end{equation*}%
and the approximation$~\widetilde{S}$ is 
\begin{eqnarray*}
\widetilde{S} &=&(I_{n_{\xi }}\otimes u_{1}^{sT})\left[ I_{n_{\xi }}\otimes
(A_{1}-\lambda _{1}^{s}I_{n_{x}})\right] ^{-1}(I_{n_{\xi }}\otimes u_{1}^{s})
\\
&\approx &I_{n_{\xi }}\otimes \left[ u_{1}^{sT}(1-\lambda
_{1}^{s})^{-1}\left( A_{1}-I_{n_{x}}\right) ^{-1}u_{1}^{s}\right] \\
&\approx &I_{n_{\xi }}\otimes \left[ u_{1}^{sT}\left( M_{1}^{s}\right)
^{-1}u_{1}^{s}\right] ,
\end{eqnarray*}%
where the second line was used in~\cite{Benner-2017-INK}. 
%[eq.~(4.28)]{Benner-2017-INK}. 
In this study, we use the third line with$~\left( M_{1}^{s}\right)
^{-1}u_{1}^{s}$ denoting an application of$~M_{1}^{s}$ to$~u_{1}^{s}$. The
ideal choice of$~u_{1}^{s}$ are the coefficients of the mean of eigenvector$%
~s$, and we consider two approximations here:\ (a) $u_{1}^{s}$\ is set as
the corresponding eigenvector of the mean matrix~$A_{1}$, or (b) $u_{1}^{s}$
is the approximation of the gPC\ coefficients of the corresponding
eigenvector updated after each step of Newton iteration (Algorithm~\ref%
{alg:line_search}). The preconditioners are thus either (a) \emph{fixed}
during Newton iteration, or (b) \emph{updated} after each step. These two
variants and our version of the mean-based preconditioner~(NMB) for problem~(%
\ref{eq:Newton-mod}) are summarized in Algorithm~\ref{alg:NMB}. Clearly, if$%
~M_{1}^{s}$ is symmetric, positive definite, so is the preconditioner~$M_{%
\text{NMB}}$, but the preconditioner loses positive definiteness if the
eigenvalue of interest is not the smallest one, cf.~(\ref{eq:M_1^s}), and
therefore, along with \textsc{MINRES}, we also use \textsc{GMRES} and
develop several preconditioners for this method.

\begin{algorithm}[hptb]
\caption{Mean-based preconditioner for the Newton iteration (NMB)}
\label{alg:NMB} 
The preconditioner $M_{\text{NMB}}:\left( \bar{r}^{(u),s},\bar{r}^{(\lambda),s}\right) \longmapsto \left( \bar{v}^{(u),s},\bar{v}^{(\lambda),s}\right) $ %for~(\ref{eq:Newton-mod}) 
is defined~as 
\begin{equation}
\left[ 
\begin{array}{cc}
I_{n_{\xi }}\otimes M^s_{1} & 0 \\ 
0 & I_{n_{\xi }}\otimes \left[ w^{s,(n)T} \left( M_{1}^s \right)^{-1}w^{s,(n)}\right] %
\end{array}%
\right] \left[ 
\begin{array}{c}
\bar{v}^{(u),s} \\ 
\bar{v}^{(\lambda),s}%
\end{array}%
\right] =\left[ 
\begin{array}{c}
\bar{r}^{(u),s} \\ 
\bar{r}^{(\lambda),s }%
\end{array}%
\right] ,  \label{eq:algNMB}
\end{equation}
where $w^{s,(n)}$ is (a)  eigenvector $w^s$ of $A_1$ corresponding to eigenvalue~$\mu^s$, 
cf.~(\ref{eq:alg-SISI-ini}), 
or (b)  the first (mean) gPC coefficients~$u^{s,(n)}_{1}$ of eigenvector~$s$ at step~$n$ of Algorithm~\ref{alg:line_search}.
\end{algorithm}

Next, we propose a variant of so-called, \emph{constraint preconditioner},
cf.~\cite{Keller-2000-CPI}, 
\begin{equation*}
\left[ 
\begin{array}{cc}
\widetilde{A} & \widetilde{B}^{T} \\ 
\widetilde{B} & 0%
\end{array}%
\right] .
\end{equation*}%
Similarly as above, both $\widetilde{A}$ and $\widetilde{B}$ are
approximations of the blocks in~(\ref{eq:jac-scheme}). The preconditioner is
clearly indefinite (which also precludes use of \textsc{MINRES}).\ Our
variant of the constraint mean-based preconditioner~(cMB) is listed as
Algorithm~\ref{alg:cMB}.

\begin{algorithm}[hptb]
\caption{Constraint mean-based preconditioner (cMB)}
\label{alg:cMB} 
The preconditioner $M_{\text{cMB}}:\left( \bar{r}^{(u),s},\bar{r}^{(\lambda),s}\right) \longmapsto \left( \bar{v}^{(u),s},\bar{v}^{(\lambda),s}\right) $\ 
is defined as 
\begin{equation}
\left[ 
\begin{array}{cc}
I_{n_{\xi }}\otimes M_{1}^s & - I_{n_{\xi }}\otimes w^{s,(n)} \\ 
- I_{n_{\xi }}\otimes w^{s,(n)T} & 0%
\end{array}%
\right] \left[ 
\begin{array}{c}
\bar{v}^{(u),s} \\ 
\bar{v}^{(\lambda),s}%
\end{array}%
\right] =\left[ 
\begin{array}{c}
\bar{r}^{(u),s} \\ 
\bar{r}^{(\lambda),s}%
\end{array}%
\right] ,  \label{eq:algcMB}
\end{equation}%
where $w^{s,(n)}$ is set as in Algorithm~\ref{alg:NMB}.
\end{algorithm}

In an analogy to Algorithm~\ref{alg:MB} and~(\ref{eq:MB-matricized}), the
action of the preconditioners from Algorithms~\ref{alg:NMB} and~\ref{alg:cMB}
can be equivalently obtained by solving 
\begin{equation}
\mathcal{M}_{1}\left[ 
\begin{array}{c}
\bar{V}^{(u),s} \\ 
\bar{V}^{(\lambda ),s}%
\end{array}%
\right] =\left[ 
\begin{array}{c}
\bar{R}^{(u),s} \\ 
\bar{R}^{(\lambda ),s}%
\end{array}%
\right] ,  \label{eq:MB-matricized-N}
\end{equation}%
where $\mathcal{M}_{1}$ is the deterministic part the preconditioners from~(%
\ref{eq:algNMB}) or~(\ref{eq:algcMB}), that is 
\begin{equation*}
\mathcal{M}_{1}=\left[ 
\begin{array}{cc}
M_{1}^{s} & 0 \\ 
0 & w^{s,(n)T}\left( M_{1}^{s}\right) ^{-1}w^{s(n)}%
\end{array}%
\right] \quad \text{or}\quad \mathcal{M}_{1}=\left[ 
\begin{array}{cc}
M_{1}^{s} & -w^{s(n)} \\ 
-w^{s,(n)T} & 0%
\end{array}%
\right] .
\end{equation*}

We also formulate a constraint version of the preconditioner from Algorithm~%
\ref{alg:hGS}, which is called a constraint hierarchical Gauss-Seidel
preconditioner~(chGS) and is formulated as Algorithm~\ref{alg:chGS}--\ref%
{alg:chGS_cont}. There are two components of the preconditioner. The first
component consists of block-diagonal solves with blocks of varying sizes
computed just as in Algorithm~\ref{alg:cMB}, resp.~(\ref{eq:MB-matricized-N}%
). The second component is used in the setup of the right-hand sides for the
solves and consists of matrix-vector products by certain subblocks of the
stochastic Jacobian matrices by vectors of corresponding sizes. 
An example of matrix-vector product with a subblock of the stochastic Jacobian matrix 
is given in Appendix~\ref{sec:mvp}. The preconditioner is formulated as
Algorithm~\ref{alg:chGS}--\ref{alg:chGS_cont}, and a scheme of the splitting
operator is illustrated by the right panel of Figure~\ref{fig:J-hierarchy}.
We also note that, since the initial guess is zero, 
%in Algorithm~\ref{alg:chGS}--\ref{alg:chGS_cont}, 
the multiplications by$~\mathcal{F}_{1}$\ and $\mathcal{F}_{d+1}$ vanish
from~(\ref{eq:algchGS1})--(\ref{eq:algchGS2}).

\begin{algorithm}[hptb]
\caption{Constraint hierarchical Gauss-Seidel preconditioner (chGS)}
\label{alg:chGS} 

The preconditioner $M_{\text{chGS}}:\left( \bar{r}^{(u),s},\bar{r}^{(\lambda
),s}\right) \longmapsto \left( \bar{v}^{(u),s},\bar{v}^{(\lambda ),s}\right) 
$\ 
is defined as follows.
\begin{algorithmic}[1]
\State Set the initial solution $\left( \bar{v}^{(u),s},\bar{v}^{(\lambda
),s}\right) $ to zero and update in the following steps:
\State Solve 
\begin{equation}
\mathcal{M}_{1}\left[ 
\begin{array}{c}
v_{1}^{(u),s} \\ 
v_{1}^{(\lambda ),s}%
\end{array}%
\right] =\left[ 
\begin{array}{c}
r_{1}^{(u),s} \\ 
r_{1}^{(\lambda ),s}%
\end{array}%
\right] -\mathcal{F}_{1}\left[ 
\begin{array}{c}
v_{\left( 2:n_{\xi }\right) }^{(u),s} \\ 
v_{\left( 2:n_{\xi }\right) }^{(\lambda ),s}%
\end{array}%
\right],   \label{eq:algchGS1}
\end{equation}%
where 
\begin{eqnarray*}
\mathcal{M}_{1} &=&\left[ 
\begin{array}{cc}
M_{1}^{s} & - w^{s,(n)} \\ 
- w^{s,(n)T} & 0%
\end{array}%
\right], \text{where } w^{s,(n)} \text{ is set as in Algorithm~\ref{alg:NMB},} \\
\mathcal{F}_{1} &=&\left[  \begin{array}{cc} \sum_{t\in \mathcal{I}_{t}}\left[ h_{t,\left(1\right)\left( 2:n_{\xi }\right) }\right]  \otimes A_{t} -  \sum_{t\in \mathcal{I}_{t}}\left[ h_{t,\left(1\right) \left( 2:n_{\xi }\right) }\right] \otimes \lambda_t^{s,(n)} I_{n_x} & \mathcal{G}_{1} \\ 
\mathcal{H}_{1} & 0%
\end{array}%
\right] ,  \\
\mathcal{G}_{1} &=&\sum_{t\in \mathcal{I}_{t}}\left[ h_{t,\left(1\right) \left( 2:n_{\xi }\right) }\right]
\otimes w_t^{s,(n)} ,\\
\mathcal{H}_{1} &=&\sum_{t\in \mathcal{I}_{t}}\left[ h_{t,\left(1\right)\left( 2:n_{\xi }\right) }\right]
\otimes (w_t^{s,(n)})^T, 
\end{eqnarray*}
where $w^{s,(n)}$ the eigenvector~$s$ at step~$n$ of Algorithm~\ref{alg:line_search}.
\For{$d=1,\ldots, p-1$}
\State Set $\ell =\left( n_{\ell }+1:n_{u}\right) ,\text{ where }n_{\ell }=\binom{n_{\xi
}+d-1}{d-1}\text{ and }n_{u}=\binom{n_{\xi }+d}{d}$.
\State  Solve 
\begin{equation}
\mathcal{M}_{d+1}\left[ 
\begin{array}{c}
v_{\left( \ell \right) }^{(u),s} \\ 
v_{\left( \ell \right) }^{(\lambda ),s}%
\end{array}%
\right] =\left[ 
\begin{array}{c}
r_{\left( \ell \right)}^{(u),s} \\ 
r_{\left( \ell \right) }^{(\lambda ),s}%
\end{array}%
\right] -\mathcal{E}_{d+1}\left[ 
\begin{array}{c}
v_{(1:n_{\ell })}^{(u),s} \\ 
v_{(1:n_{\ell })}^{(\lambda ),s}%
\end{array}%
\right] -\mathcal{F}_{d+1}\left[ 
\begin{array}{c}
v_{(n_{u}+1:n_{\xi })}^{(u),s} \\ 
v_{(n_{u}+1:n_{\xi })}^{(\lambda ),s}%
\end{array}%
\right] ,  \label{eq:algchGS2}
\end{equation}

where 
\begin{eqnarray*}
\mathcal{M}_{d+1} &=&\left( I_{n_{u}-n_{\ell }} \!\otimes\! \left[ 
\begin{array}{cc}
M_{1}^s\!\! & \!\!-w^{s,(n)} \!\!\\ 
-w^{s,(n)T}\!\! & 0%
\end{array}%
\right] \right) , \text{where } w^{s,(n)} \text{ is set as in Algorithm~\ref{alg:NMB},}  \\
\mathcal{E}_{d+1} &=&\left[  \begin{array}{cc} \sum_{t\in \mathcal{I}_{t}} \left[ h_{t,(\ell )\left( 1:n_{\ell }\right) }
\right]  \otimes A_{t} -  \sum_{t\in \mathcal{I}_{t}} \left[ h_{t,(\ell )\left( 1:n_{\ell }\right) }
\right]  \otimes \lambda_t^{s,(n)} I_{n_x} & \mathcal{G}_{d+1}^{\mathcal{E}} \\ 
\mathcal{H}_{d+1}^{\mathcal{E}} & 0%
\end{array}%
\right] , \\
\mathcal{G}_{d+1}^{\mathcal{E}} &=&\sum_{t\in \mathcal{I}_{t}} \left[ h_{t,(\ell )\left( 1:n_{\ell }\right) }
\right]  \otimes w_t^{s,(n)} , \\
\mathcal{H}_{d+1}^{\mathcal{E}} &=&\sum_{t\in \mathcal{I}_{t}} \left[ h_{t,(\ell )\left( 1:n_{\ell }\right) }
\right] 
\otimes (w_t^{s,(n)})^T , \\
\mathcal{F}_{d+1} &=&\left[  \!\begin{array}{cc} \sum_{t\in \mathcal{I}_{t}}\! \left[ h_{t,\left( \ell
\right) \left( n_{u}+1:n_{\xi }\right) }\right] \! \otimes \!A_{t}\! -\!  \sum_{t\in \mathcal{I}_{t}}\! \left[ h_{t,\left( \ell
\right) \left( n_{u}+1:n_{\xi }\right) }\right] \! \otimes\! \lambda_t^{s,(n)} I_{n_x} & \!\mathcal{G}_{d+1}^{\mathcal{F}}\! \\ 
\mathcal{H}_{d+1}^{\mathcal{F}} & 0%
\end{array}%
\right] , \\
\mathcal{G}_{d+1}^{\mathcal{F}} &=&\sum_{t\in \mathcal{I}_{t}}\left[ h_{t,\left( \ell
\right) \left( n_{u}+1:n_{\xi }\right) }\right] \otimes w_t^{s,(n)} , \\
\mathcal{H}_{d+1}^{\mathcal{F}} &=&\sum_{t\in \mathcal{I}_{t}} \left[ h_{t,\left( \ell
\right) \left( n_{u}+1:n_{\xi }\right) }\right]
\otimes (w_t^{s,(n)})^T .
\end{eqnarray*}
\EndFor
\algstore{split_here}
\end{algorithmic}
\end{algorithm}

\begin{algorithm}[hptb]
\caption{Constraint hierarchical Gauss-Seidel preconditioner (chGS), continued}
\label{alg:chGS_cont} 
\begin{algorithmic}[1]
\algrestore{split_here}
\State Set $\ell =\left( n_{u}+1:n_{\xi }\right)$.
\State Solve%
\begin{eqnarray*}
\mathcal{M}_{p+1}\left[ 
\begin{array}{c}
v_{\left( \ell \right) }^{(u),s} \\ 
v_{\left( \ell \right) }^{(\lambda ),s}%
\end{array}%
\right] &=&\left[ 
\begin{array}{c}
r_{\left( \ell \right) }^{(u),s} \\ 
r_{\left( \ell \right) }^{(\lambda ),s}%
\end{array}%
\right] -\mathcal{E}_{p+1}\left[ 
\begin{array}{c}
v_{\left( 1:n_{u}\right) }^{(u),s} \\ 
v_{\left( 1:n_{u}\right) }^{(\lambda ),s}%
\end{array}%
\right] , \\
\end{eqnarray*}  %
where 
\begin{eqnarray*}
\mathcal{M}_{p+1} &=&\left( I_{n_{\xi }-n_{u}}\!\otimes\!\! \left[ 
\begin{array}{cc}
\!M_{1}^s \!& \!\!-w^{s,(n)}\!\! \\ 
\!-w^{s,(n)T} \!& 0\!%
\end{array}%
\right] \right) , \text{where } w^{s,(n)} \text{ is set as in Algorithm~\ref{alg:NMB},}  \\
\mathcal{E}_{p+1} &=&\left[  \begin{array}{cc} \sum_{t\in \mathcal{I}_{t}} \left[ h_{t,(\ell )\left( 1:n_u \right) }
\right]  \otimes A_{t} -  \sum_{t\in \mathcal{I}_{t}} \left[ h_{t,(\ell )\left( 1:n_u \right) }
\right]  \otimes \lambda_t^{s,(n)} I_{n_x} & \mathcal{G}_{p+1}^{\mathcal{E}} \\ 
\mathcal{H}_{p+1}^{\mathcal{E}} & 0%
\end{array}%
\right] ,  \\
\mathcal{G}_{p+1}^{\mathcal{E}} &=&\sum_{t\in \mathcal{I}_{t}} \left[ h_{t,(\ell )\left( 1:n_u\right) }
\right]  \otimes w_t^{s,(n)},\\  %\left( \left[ h_{(2:n_{\xi }),11}\right] \otimes
%I_{n_{x}}\right) \bar{u}_{1}^{s,(n)} \\
\mathcal{H}_{p+1}^{\mathcal{E}} &=&\sum_{t\in \mathcal{I}_{t}} \left[ h_{t,(\ell )\left( 1:n_u \right) }
\right] 
\otimes (w_t^{s,(n)})^T .
\end{eqnarray*}
\For{$d=p-1,\ldots, 1$}
\State Set $\ell =\left( n_{\ell }+1:n_{u}\right) ,\text{ where }n_{\ell }=\binom{n_{\xi
}+d-1}{d-1}\text{ and }n_{u}=\binom{n_{\xi }+d}{d}$.
\State  Solve~(\ref{eq:algchGS2}).
\EndFor
\State Solve~(\ref{eq:algchGS1}).
\end{algorithmic}
\end{algorithm}

\section{Numerical experiments}
We implemented the methods in \textsc{Matlab}, and in this section we present the results of numerical experiments in which the proposed inexact solvers are applied to two benchmark problems: a diffusion problem with stochastic coefficient and stiffness of Mindlin plate with stochastic Young's modulus.

\subsection{Stochastic diffusion problem with lognormal coefficient} 
\label{sec:numerical} For the first benchmark problem we consider
the elliptic equation with stochastic coefficient and deterministic Dirichlet boundary condition 
\begin{eqnarray*}
-\nabla \cdot \left( a(x,\xi )\nabla u(x,\xi )\right) &=&\lambda (\xi
)u(x,\xi )\quad \text{in }D\times \Gamma , \\
u(x,\xi ) &=&0\qquad \qquad \quad \;\text{on }\partial D\times \Gamma ,
\end{eqnarray*}%
where $D$ is a two-dimensional physical domain. The uncertainty in the model
is introduced by the stochastic expansion of the diffusion coefficient,
considered as 
\begin{equation}
a(x,\xi )=\sum_{\ell =1}^{n_{a}}a_{\ell }(x)\psi _{\ell }(\xi ),
\label{eq:a-stoch_exp}
\end{equation}%
to be a truncated lognormal process transformed from the underlying Gaussian
process~\cite{Ghanem-1999-NGS}. That it, $\psi _{\ell }(\xi )$, $\ell
=1,\dots ,n_{a}$, is a set of Hermite polynomials and, denoting the
coefficients of the Karhunen-Lo\`{e}ve expansion of the Gaussian process by$%
~g_{j}(x)$ and $\eta _{j}=\xi _{j}-g_{j}$, $j=1,\dots ,m_{\xi }$, the
coefficients in expansion~(\ref{eq:a-stoch_exp}) are computed as 
\begin{equation*}
a_{\ell }(x)=\frac{\mathbb{E}\left[ \psi _{\ell }(\eta )\right] }{\mathbb{E}%
\left[ \psi _{\ell }^{2}(\eta )\right] }\exp \left[ g_{0}+\frac{1}{2}%
\sum_{j=1}^{m_{\xi }}\left( g_{j}(x)\right) ^{2}\right] .
\end{equation*}%
%The covariance function of the Gaussian field was chosen to be 
%\begin{equation}
%C\left( x_{1},x_{2}\right) =\sigma _{g}^{2}\exp \left( -\frac{\left\Vert
%x_{1}-x_{2}\right\Vert _{2}}{L_{\text{corr}}}\right) , \label{eq:covariance}
%\end{equation}%
%where$~L_{\text{corr}}$ is the correlation length and $\sigma _{g}$ is the standard deviation of the Gaussian random field. 
The covariance function of the Gaussian field, for points $X_{1}=(x_{1}%
,y_{1})$ and $X_{2}=(x_{2},y_{2})$ in$~D$, was chosen to be
\begin{equation}
C\left(  X_{1},X_{2}\right)  =\sigma_{g}^{2}\exp\left(  -\frac{\left\vert
x_{2}-x_{1}\right\vert }{L_{x}}-\frac{\left\vert y_{2}-y_{1}\right\vert
}{L_{y}}\right)  , \label{eq:covariance}%
\end{equation}
where$~L_{x}$ and $L_{y}$ are the correlation lengths of the random variables
$\xi_{i}$, $i=1,\dots,m_{\xi}$, in the $x$ and $y$\ directions, respectively,
and $\sigma_{g}$ is the standard deviation of the Gaussian random field.
According to~\cite%
{Matthies-2005-GML}, in order to guarantee a complete representation of the
lognormal process by~(\ref{eq:a-stoch_exp}), the degree of polynomial
expansion of~$a(x,\xi )$ should be twice the degree of the expansion of the
solution. We follow the same strategy here. Therefore, the values of $n_{\xi
}$ and $n_{a}$ are, cf., e.g.~\cite[p.~87]{Ghanem-1991-SFE} or~\cite[%
Section~5.2]{Xiu-2010-NMS}, 
%\begin{equation*}
$n_{\xi }\!=\!\frac{\left( m_{\xi }+p\right) !}{m_{\xi }!p!}$, $n_{a}\!=\!\frac{%
\left( m_{\xi }+2p\right) !}{m_{\xi }!\left( 2p\right) !}$.
%\end{equation*}%
In the numerical experiments, the lognormal diffusion coefficient (\ref%
{eq:a-stoch_exp}) is parameterized using $m_{\xi }=3$ random variables. The
correlation length is $L_{\text{corr}}=2$, and the coefficient of variation $%
CoV$ of the lognormal process is set either to $0.1$ ($10\%)$ or $0.25$ ($%
25\%)$, where $CoV=\sigma /a_{1}$, the ratio of the standard deviation~$%
\sigma $ and the mean of the diffusion coefficient$~a_{1}$. For the gPC
expansion of eigenvalues/eigenvectors~(\ref{eq:sol_mat}), the maximal degree
of gPC\ expansion is $p=3$, so then $n_{\xi }=20$ and $n_{a}=84$.

Finite element discretization leads to a generalized eigenvalue problem 
\begin{equation}
K(\xi )u=\lambda Mu,  \label{eq:gen-eig}
\end{equation}%
where $K(\xi )=\sum_{\ell =1}^{n_{a}}K_{\ell }\psi _{\ell }(\xi )$ is the
stochastic expansion of the stiffness matrix, and the mass matrix$~M$ is
deterministic. Using Cholesky factorization $M=LL^{T}$, the generalized
eigenvalue problem~(\ref{eq:gen-eig}) can be transformed into the standard
form 
\begin{equation}
A(\xi )w=\lambda w,  \label{eq:diff_eig}
\end{equation}%
where $u=L^{-T}w$ and the expansion of $A$ corresponding to~(\ref%
{eq:stoch-exp-A}) is 
\begin{equation}
A=\sum_{\ell =1}^{n_{a}}A_{\ell }\psi _{\ell }(\xi )=\sum_{\ell =1}^{n_{a}} 
\left[ L^{-1}K_{\ell }L^{-T}\right] \psi _{\ell }(\xi ).
\label{eq:diff_coeff}
\end{equation}

\begin{figure}[tbh]
\centering
\includegraphics[angle=0, scale=0.38]{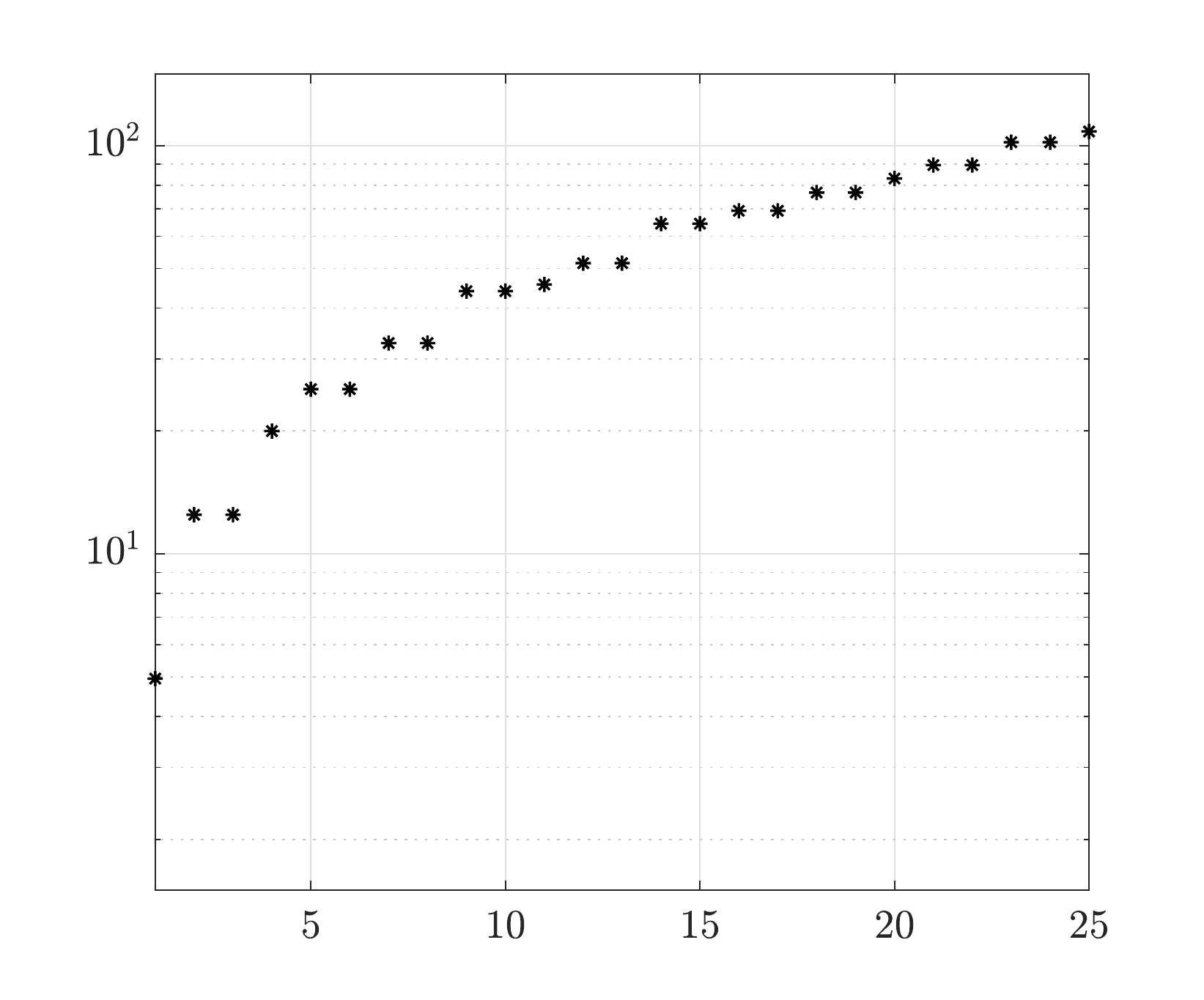}
\caption{The smallest 25 eigenvalues of the mean matrix $A_{1}$.}
\label{fig:eigval_mean}
\end{figure}

We consider the physical domain $D=[-1,1]^{2},$ discretized using a
structured grid using $256$ bilinear finite elements, that is with $225$
nodes interior to$~D$, which determines the size of matrices$~A_{\ell }$ in~(%
\ref{eq:diff_coeff}). The$~25$ smallest eigenvalues of the mean matrix $%
A_{1} $ are displayed\ in Figure~\ref{fig:eigval_mean}. For the quadrature
rule, in Section~\ref{sec:sampling}, we use Smolyak sparse grid with
Gauss-Hermite quadrature and grid level~$4$, and $10^{4}$ samples for the
Monte Carlo method. With these settings, the size of~$h_{\ell ,kj}$ in~(\ref%
{eq:cijk}) was $84\times 20\times 20$ with $806$ nonzeros, and there were $%
69 $ points on the sparse grid.

\subparagraph{Inexact stochastic inverse subspace iteration}

\begin{figure}[tbh]
\vspace{-6mm}
\centering
\subfloat{	\includegraphics[angle=0, scale=0.582]{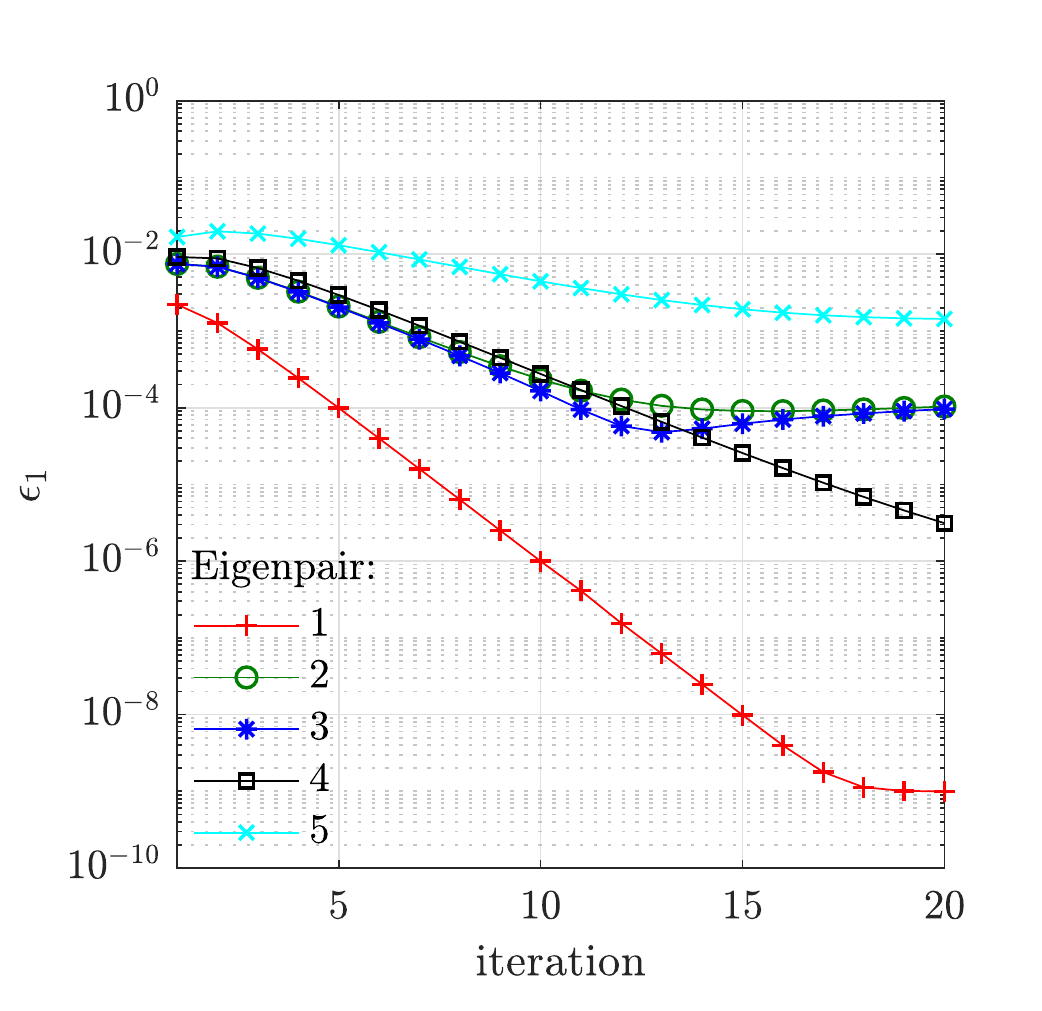}
	\label{fig:sisi_decay_a}
	} 
\subfloat{	\includegraphics[angle=0, scale=0.582]{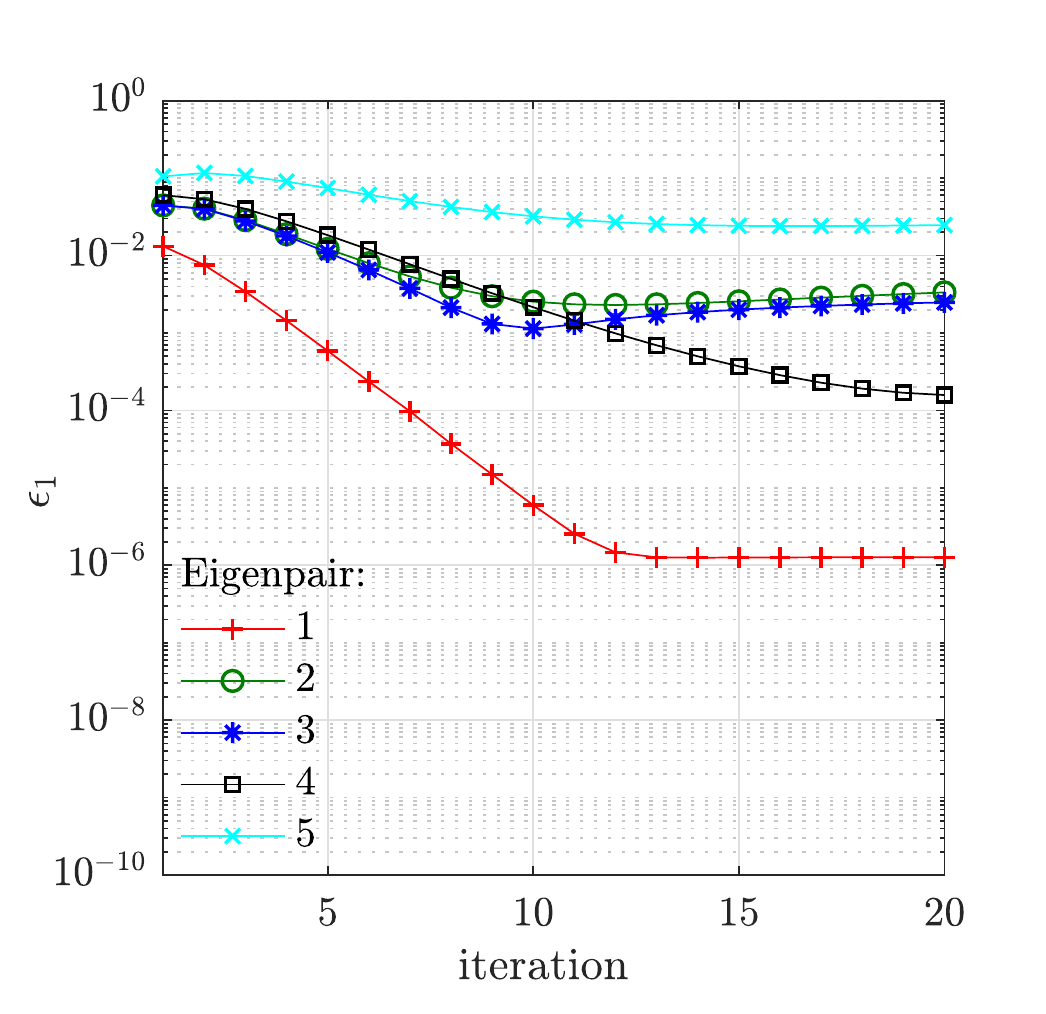}
	\label{fig:sisi_decay_b}
	}\vspace{-5.5mm}\newline
\subfloat{	\includegraphics[angle=0, scale=0.582]{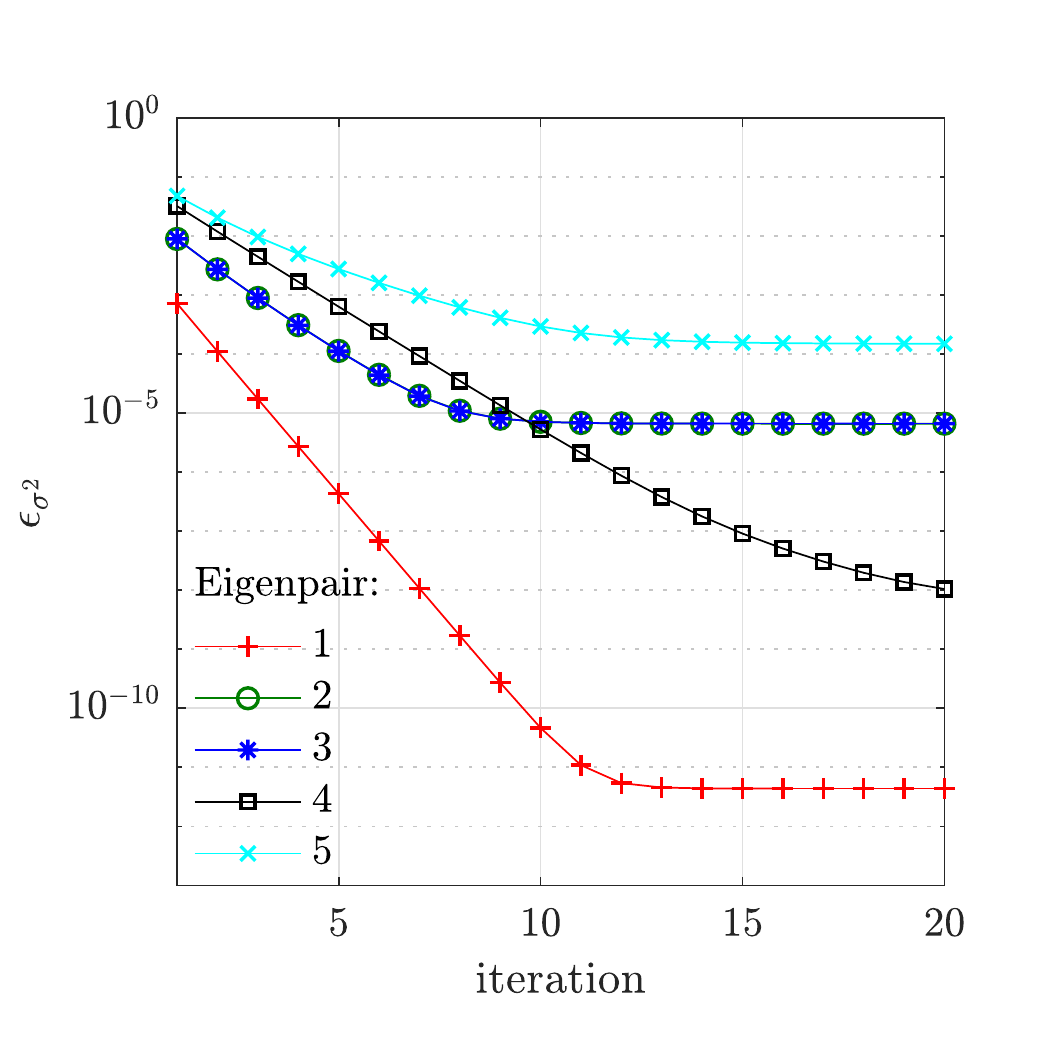}
	\label{fig:sisi_decay_c}
	} 
\subfloat{	\includegraphics[angle=0, scale=0.582]{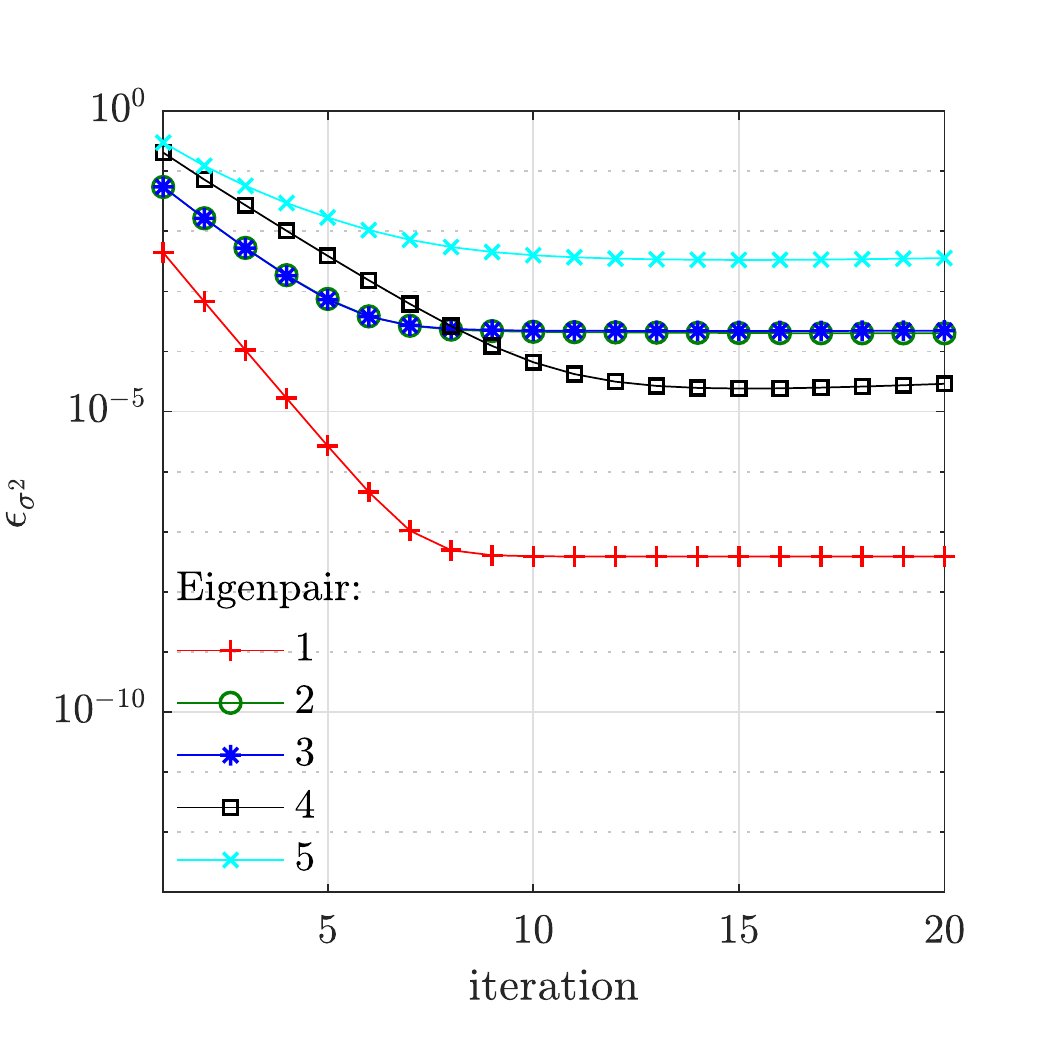}
	\label{fig:sisi_decay_d}
	}	\vspace{-4mm}
\caption{Convergence history of the inexact stochastic inverse subspace
iteration in terms of indicators $\protect\epsilon_1$ (top) and $\protect%
\epsilon_{\protect\sigma^2}$ (bottom) defined by~(\protect\ref{eq:eps}) with 
$CoV=10\%$ (left) and $25\%$ (right).}
\label{fig:sisi_decay}
\end{figure}

First, we examine the performance of the inexact stochastic inverse subspace
iteration (SISI) from Algorithm~\ref{alg:isisi} for computing the five
smallest eigenvalues and corresponding eigenvectors of problem~(\ref%
{eq:diff_eig}). Linear systems~(\ref{eq:alg-SISI-solve}) are solved using
the PCG method with the mean-based preconditioner (Algorithm~\ref{alg:MB})
and the hierarchical Gauss-Seidel preconditioner (Algorithm~\ref{alg:hGS}).
We ran the SISI algorithm with a fixed number of steps set to $20$. Figure~%
\ref{fig:sisi_decay} illustrates convergence history in terms of the two
error indicators$~\epsilon _{1}$ and$~\epsilon _{\sigma ^{2}}$ from$~$(\ref%
{eq:eps}) with $CoV=10\%$ (left panels) and $25\%$ (right panels). The plots
were generated using the hGS preconditioner with $p_{t}=2$ ($\mathcal{I}%
_{t}=\{1,\ldots ,10\}$), but convergence with other preconditioners was
virtually identical.

Next, we examine performance of PCG with the two preconditioners used to
solve linear systems~(\ref{eq:alg-SISI-solve}) with zero initial guess and
stopping criterion~(\ref{eq:pcg-stop}). We computed the five smallest
eigenvalues using 20 steps of the inexact SISI method. Table~\ref%
{tab:pcg_iters} shows the number of the PCG iterations required by the
inexact solves, averaged over the 20 steps of the inexact SISI method for
the model eigenvalue problem with $CoV=10\%$ and$~25\%$. Specifically, we
compare the mean-based preconditioner from Algorithm~\ref{alg:MB} and the
hGS preconditioner from Algorithm~\ref{alg:hGS} with varying level of
truncation of the matrix-vector multiplications ($p_{t}=\{0,1,2\}$ and $%
p_{t}=3$, i.e., no truncation). In both preconditioners we used Cholesky
factorization of $A_{1}$ for the solves with$~M_{1}$. We note that with $%
p_{t}=0$\ the hGS preconditioner reduces to the mean-based preconditioner.
In both cases $CoV=10\%$ and $25\%$ the hGS preconditioner outperforms the
mean-based preconditioner in terms of the number of PCG iterations for each
of the five eigenpairs. Table~\ref{tab:pcg_iters} also shows that solving
the eigenvalue problem with higher $CoV$ leads to only a slight increase in
the number of iterations.

\begin{table}[!h]
\caption{Average number of PCG iterations for computing the five smallest
eigenvalues and corresponding eigenvectors of the diffusion problem with $%
CoV =10\%$ (left) and $25\%$ (right) using inexact stochastic inverse
subspace iteration (Algorithm~\protect\ref{alg:isisi}).}
\label{tab:pcg_iters}
\begin{center}
{\footnotesize \renewcommand{\arraystretch}{1.3} 
\begin{tabular}{|c|ccccc|ccccc|}
\hline
& \multicolumn{5}{c}{$CoV=10\%$} & \multicolumn{5}{|c|}{$CoV=25\%$} \\ 
\hline
Preconditioner & 1st & 2nd & 3rd & 4th & 5th & 1st & 2nd & 3rd & 4th & 5th
\\ \hline
MB & 6.45 & 3.90 & 3.90 & 4.60 & 3.75 & 8.60 & 5.55 & 5.55 & 6.05 & 4.75 \\ 
hGS ($p_t=1$) & 3.10 & 1.95 & 1.95 & 2.25 & 1.95 & 3.65 & 2.75 & 2.75 & 2.65
& 2.00 \\ 
hGS ($p_t=2$) & 2.35 & 1.70 & 1.70 & 1.65 & 1.00 & 2.60 & 1.90 & 1.90 & 1.85
& 1.75 \\ 
hGS (no trunc.) & 2.15 & 1.00 & 1.00 & 1.45 & 1.00 & 2.60 & 1.80 & 1.80 & 
1.75 & 1.65 \\ \hline
\end{tabular}
}
\end{center}
\end{table}

\subparagraph{Newton iteration}

\begin{figure}[!t]
\vspace{-5mm}
\centering
\subfloat{	\includegraphics[angle=0, scale=0.60]{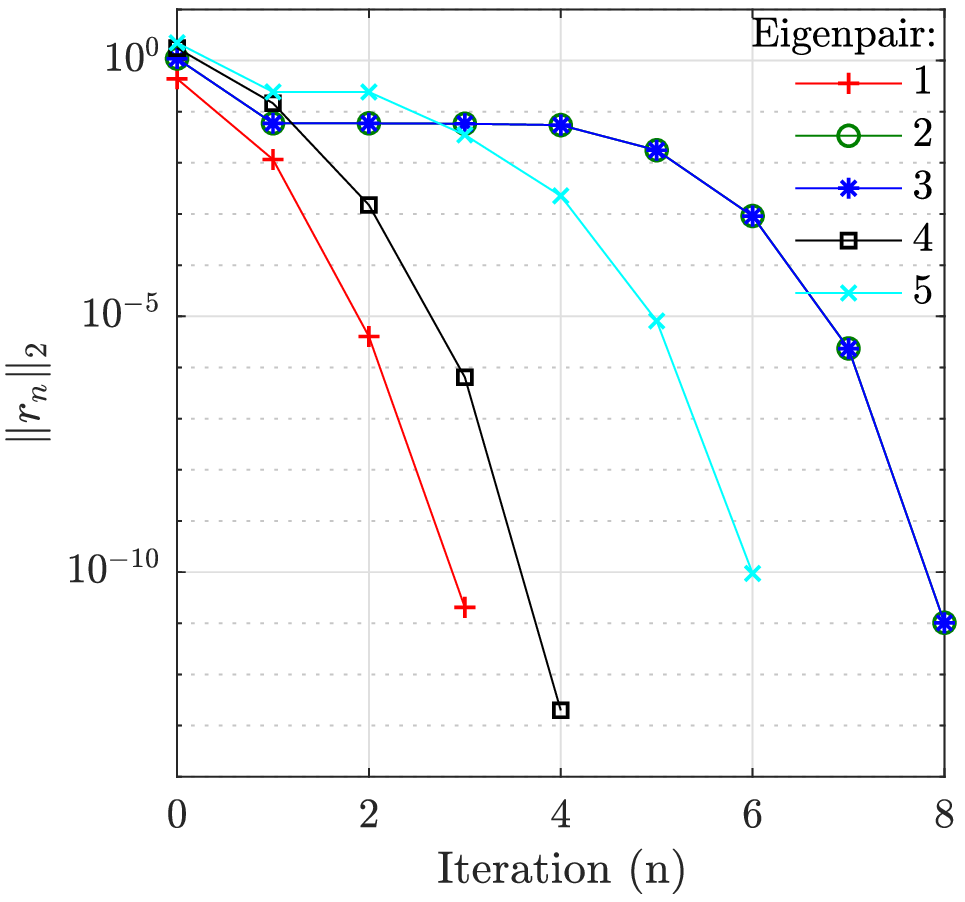}
	\label{fig:nlres_decay_a}
	} 
\subfloat{	\includegraphics[angle=0, scale=0.60]{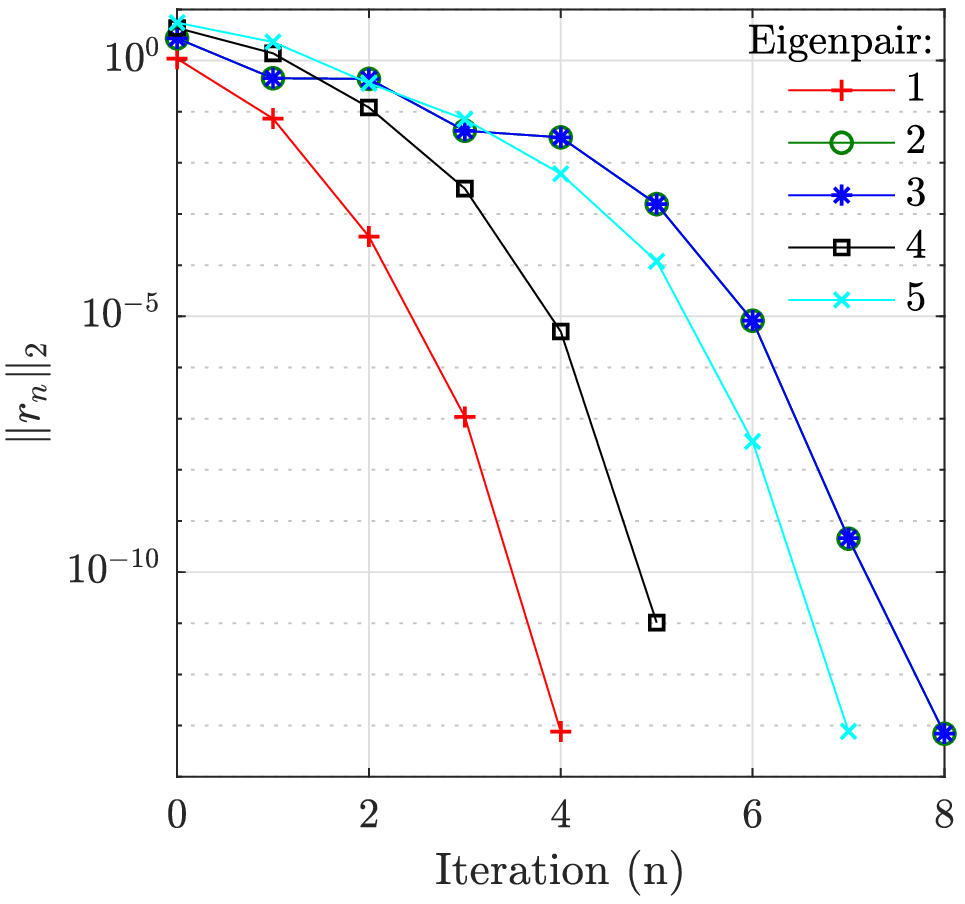}
	\label{fig:nlres_decay_b}
	}
	     \vspace{-5mm}
\caption{Convergence history in terms of the nonlinear residual $\| r_{n}
\|_{2}$ of the inexact line-search Newton method with $CoV=10\%$ (left) and $%
25\%$ (right).}
\label{fig:nlres_decay}
\end{figure}

Next, we examine the inexact line-search Newton method from Algorithm~\ref%
{alg:line_search} for computing the five smallest eigenvalues and
corresponding eigenvectors of problem~(\ref{eq:diff_eig}). For the
line-search method, we set $\rho =0.9$ for the backtracking and limit the
maximum number of backtracks to$~25$, and $c=0.05$. The initial guess for
the nonlinear iteration is set using the (five smallest) eigenvalues and
corresponding eigenvectors of the eigenvalue problem associated with the
mean matrix $A_{1}$ as discussed in Section~\ref{sec:line_search}. The
nonlinear iteration terminates when the norm of the residual $\Vert
r_{n}\Vert _{2}<10^{-10}$. The linear systems in Line~\ref{ln:lin_sys} in
Algorithm~\ref{alg:line_search} are solved using either MINRES or GMRES with
the mean-based preconditioner (Algorithm~\ref{alg:NMB}), constraint
mean-based preconditioner (Algorithm~\ref{alg:cMB}) and the contraint
hierarchical Gauss-Seidel preconditioner (Algorithm~\ref{alg:chGS}--\ref%
{alg:chGS_cont}). Figure~\ref{fig:nlres_decay} illustrates convergence
history of the inexact line-search Newton method in terms of norm of the
residual$~\Vert r_{n}\Vert _{2}$ with $CoV=10\%$ (left panel) and$~25\%$%
~(right panel). The plots were generated using GMRES with the chGS
preconditioner (Algorithm~\ref{alg:chGS}--\ref{alg:chGS_cont}) with $p_{t}=2$
($\mathcal{I}_{t}=\{1,\ldots ,10\}$), but convergence with other
preconditioners was virtually identical.

\begin{table}[!t]
\caption{Average number of MINRES/GMRES iterations for computing the five
smallest eigenvalues and corresponding eigenvectors of the diffusion problem
with $CoV = 10\%$ (left) and $25\%$ (right) using inexact line-search Newton
method (Algorithm~\protect\ref{alg:line_search}) with the stopping criteria $\Vert
r_{n}\Vert _{2}<10^{-10}$.}
\label{tab:nl_iters}
\begin{center}
{\footnotesize \renewcommand{\arraystretch}{1.3} 
\begin{tabular}{|c|rrrrr|rrrrr|}
\hline
& \multicolumn{5}{c|}{$CoV=10\%$} & \multicolumn{5}{c|}{$CoV=25\%$} \\ 
\hline
{\scriptsize {Preconditioner}} & \multicolumn{1}{c}{1st} & \multicolumn{1}{c}{2nd} & \multicolumn{1}{c}{3rd} & \multicolumn{1}{c}{4th} & \multicolumn{1}{c}{5th} & \multicolumn{1}{|c}{1st} & \multicolumn{1}{c}{2nd} & 
\multicolumn{1}{c}{3rd} & \multicolumn{1}{c}{4th} & \multicolumn{1}{c|}{5th} \\ \hline
NMB{\tiny {\ (MINRES)}} & 11.5 & 59.3 & 60.2 & 23.3 & 217.6 & 13.3 & 110.0 & 
109.4 & 49.3 & 142.9 \\ 
NMB{\tiny {\ (fixed)}} & 11.3 & 71.5 & 59.9 & 29.6 & 120.5 & 15.2 & 79.3 & 
79.5 & 43.8 & 101.1 \\ 
NMB{\tiny {\ (updated)}} & 13.3 & 28.9 & 27.8 & 16.2 & 43.0 & 19.0 & 68.9 & 
64.5 & 87.3 & 122.9 \\ 
cMB{\tiny {\ (fixed)}} & 7.0 & 37.9 & 39.5 & 8.8 & 28.1 & 13.3 & 56.6 & 
56.6 & 14.6 & 32.4 \\ 
cMB{\tiny {\ (updated)}} & 4.3 & 24.7 & 25.4 & 5.3 & 28.0 & 7.8 & 33.4 & 
33.1 & 8.6 & 15.6 \\ 
chGS($p_t=1$) & 2.3 & 17.9 & 17.1 & 2.8 & 15.4 & 3.3 & 18.3 & 18.1 & 2.8
& 18.9 \\ 
chGS($p_t=2$) & 2.0 & 12.4 & 12.5 & 2.0 & 8.5 & 3.3 & 18.9 & 19.4 & 2.4
& 10.3 \\ 
chGS(full) & 2.0 & 13.8 & 13.5 & 2.0 & 12.3 & 3.3 & 15.1 & 15.1 & 2.8 & 
14.4 \\ \hline
\end{tabular}
}
\end{center}
\end{table}

Next, we compare performance of MINRES and GMRES with the preconditioners
from Algorithms~\ref{alg:NMB}--\ref{alg:chGS_cont} used to solve linear
systems at Line~\ref{ln:lin_sys} in Algorithm~\ref{alg:line_search} with
zero initial guess and the stopping criterion~(\ref{eq:gmres-stop}). Table~%
\ref{tab:nl_iters} shows the numbers of MINRES or GMRES iterations required
by the inexact solves, averaged over the number of the nonlinear steps.
Specifically, we compare the mean-based preconditioner (NMB) from Algorithm~%
\ref{alg:NMB}, contraint mean-based preconditioner (cMB) from Algorithm~\ref%
{alg:cMB} and the constraint hierarchical Gauss-Seidel preconditioner (chGS)
from Algorithm~\ref{alg:chGS}--\ref{alg:chGS_cont}. For all preconditioners,
we need to select the vector$~w^{s,(n)}$ as discussed in Algorithm~\ref%
{alg:NMB}. Choice (a) is referred to as \emph{fixed} because the vector$%
~w^{s,(n)}$ is the corresponding eigenvector of the mean matrix~$A_{1}$, and
choice (b) is referred to as \emph{updated} because 
%the approximation of mean of the corresponding eigenvector
the vector is updated after each step of Newton iteration. Only the variant
(b) was used for the chGS preconditioner. We also need to specify (the
solves with) the matrix~$M_{1}^{s}$, in particular the choice of$~\epsilon
_{M}$ in~(\ref{eq:M_1^s}). We report values of$~\epsilon _{M}$ that, in our
experience, worked best. For (both fixed) NMB and cMB, we set $\epsilon
_{M}=0.95$. For {(updated)} cMB and chGS, we set $\epsilon _{M}=1$ and use
the SVD decomposition as $\mathcal{M}_{1}=\sum_{i=1}^{\text{rank}(\mathcal{M}%
_{1})}d_{i}y_{i}z_{i}^{T}$ to solve linear systems in~(\ref%
{eq:MB-matricized-N}). If$~\mathcal{M}_{1}$ appears to be numerically
singular, the action of the inverse of$~\mathcal{M}_{1}$ is replaced by a 
\textit{pseudoinverse} $\sum_{i=1}^{\text{rank}(\mathcal{M}%
_{1})}d_{i}^{-1}z_{i}y_{i}^{T}$. We note that with $p_{t}=0$ the chGS
preconditioner reduces to the {(updated)} cMB{\ preconditioner.} With all
preconditioners the convergence was faster for simple eigenvalues, and the
iteration counts increased in the course of Newton iteration. In both cases
with $CoV=10\%$ and $25\%$ the constraint preconditioners outperform the
mean-based preconditioners, and updating the vector$~w^{s,(n)}$ 
%after each step of nonlinear iteration 
improves the convergence. The lowest iteration counts were obtained with the
chGS preconditioner, in particular with $p_{t}=2$ and full, and we
note that the computational cost with $p_{t}=2$ is lower due to the
truncation of the matrix-vector products. For these two preconditioners,
Tables~\ref{tab:nl_iters} and~\ref{tab:gmres_iters} show that solving
the eigenvalue problem with higher $CoV$ leads to only a slight increase in
the number of iterations, and for simple eigenvalues the average iteration
counts are only slightly larger than those of~SISI.

\begin{table}[!t]
\caption{The number of GMRES iterations for computing the five smallest
eigenvalues and corresponding eigenvectors of the diffusion problem with $%
CoV = 10\%$ (left) and $25\%$ (right) using inexact line-search Newton
method (Algorithm~\protect\ref{alg:line_search}) with preconditioners cMB
(top) and chGS($p_t=2$) (bottom), and with the stopping criteria $\Vert
r_{n}\Vert _{2}<10^{-10}$.}
\label{tab:gmres_iters}
\begin{center}
{\footnotesize \renewcommand{\arraystretch}{1.3} 
\begin{tabular}{|c|ccccc|ccccc|}
\hline
& \multicolumn{5}{c|}{$CoV=10\%$} & \multicolumn{5}{c|}{$CoV=25\%$} \\ 
\hline
& 1st & 2nd & 3rd & 4th & 5th & 1st & 2nd & 3rd & 4th & 5th \\ \hline
Nonlinear step & \multicolumn{10}{c|}{cMB{\tiny {\ (updated)}}} \\ \hline
1 & 2 & 2 & 2 & 1 & 1 & 2 & 2 & 2 & 1 & 1 \\ 
2 & 4 & 8 & 8 & 3 & 35 & 5 & 8 & 8 & 3 & 3 \\ 
3 & 7 & 10 & 10 & 6 & 35 & 9 & 11 & 11 & 6 & 10 \\ 
4 &  & 13 & 14 & 11 & 21 & 15 & 18 & 17 & 12 & 12 \\ 
5 &  & 23 & 26 &  & 17 &  & 34 & 34 & 21 & 16 \\ 
6 &  & 45 & 46 &  & 23 &  & 75 & 74 &  & 22 \\ 
7 &  & 72 & 72 &  & 41 &  & 86 & 86 &  & 45 \\ 
8 &  &  &  &  & 51 &  &  &  &  &  \\ \hline
Nonlinear step & \multicolumn{10}{c|}{chGS($p_t=2$)} \\ \hline
1 & 1 & 1 & 1 & 1 & 1 & 1 & 1 & 1 & 1 & 1 \\ 
2 & 2 & 5 & 5 & 1 & 5 & 2 & 6 & 6 & 1 & 2 \\ 
3 & 3 & 6 & 6 & 2 & 5 & 4 & 4 & 4 & 2 & 5 \\ 
4 &  & 6 & 6 & 4 & 8 & 6 & 10 & 10 & 3 & 7 \\ 
5 &  & 7 & 7 &  & 10 &  & 12 & 12 & 5 & 11 \\ 
6 &  & 12 & 12 &  & 22 &  & 18 & 22 &  & 14 \\ 
7 &  & 21 & 21 &  &  &  & 45 & 45 &  & 32 \\ 
8 &  & 41 & 42 &  &  &  & 55 & 55 &  &  \\ \hline
\end{tabular}
}
\end{center}
\end{table}

\begin{figure}[!h]
\vspace{-5mm}
\centering
\subfloat{	\includegraphics[angle=0, scale=0.45]{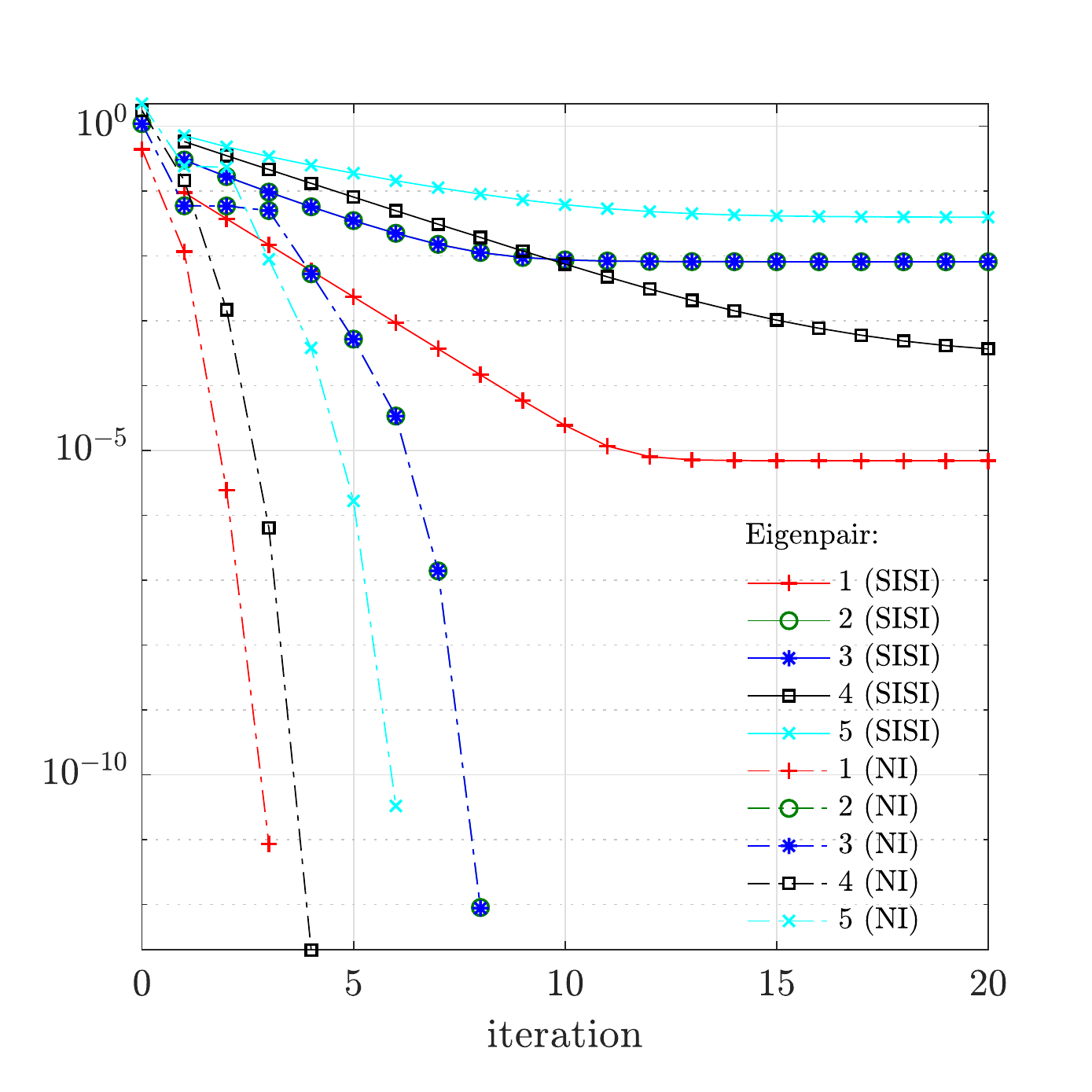}
\label{fig:sisi-ni_comparison_a}
	} 
\subfloat{	\includegraphics[angle=0, scale=0.45]{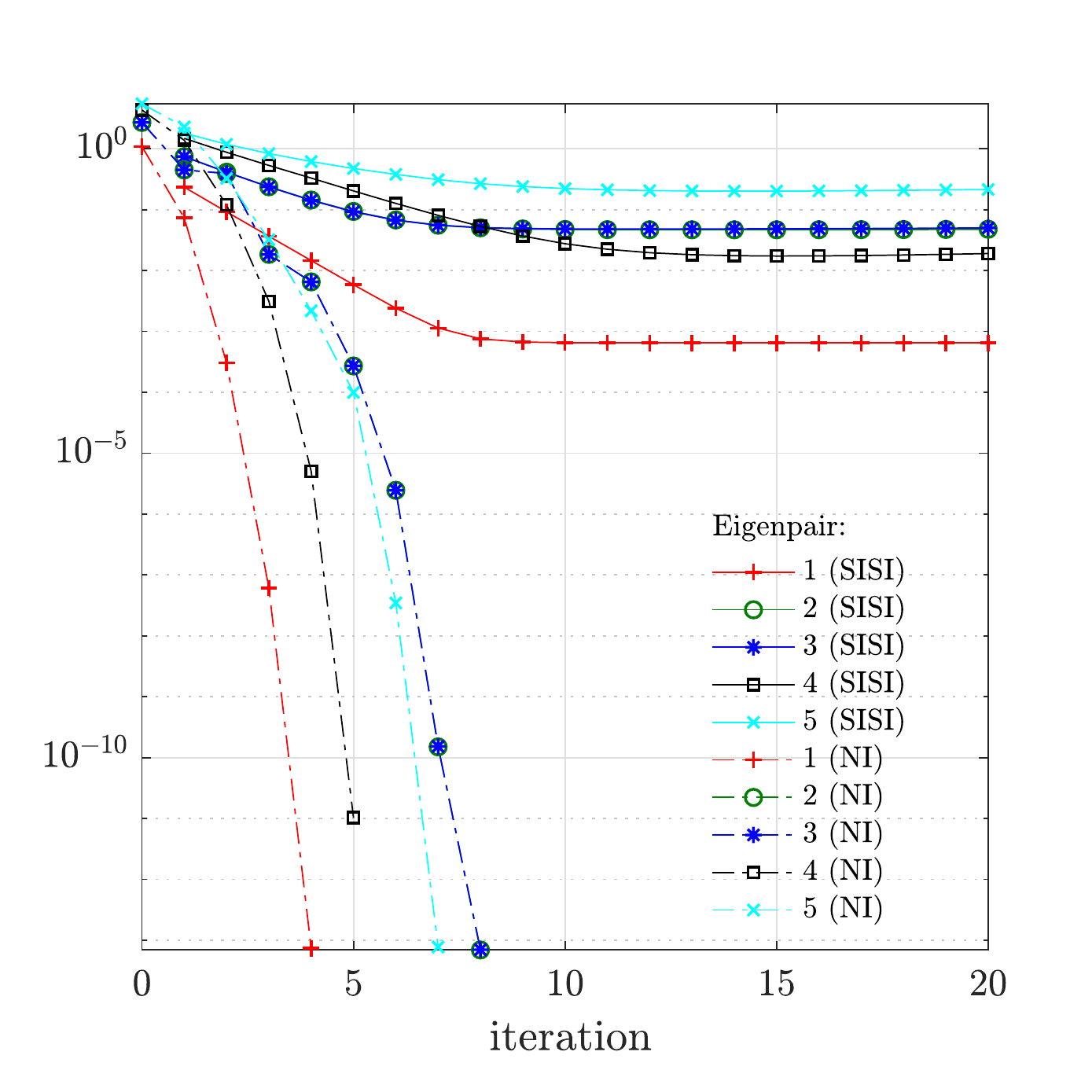}
\label{fig:sisi-ni_comparison_b}
	}
	     \vspace{-5mm}
\caption{Comparison of convergence of the inexact stochastic inverse
subspace iteration~(SISI) in terms of residual indicator~$\| \bar{\protect%
\widetilde{r}}^{s,\left( n\right) } \|_{2}$ and the inexact line-search
Newton method~(NI) in terms of~$\| F(\bar{u}^{s,(n)},\bar{\protect\lambda}%
^{s,(n)}) \|_{2}$ with $CoV=10\%$ (left) and $25\%$ (right).}
\label{fig:sisi-ni_comparison}
\end{figure}

\begin{table}[!t]
\caption{The first $10$ coefficients of the gPC expansion of the smallest
eigenvalue of the diffusion problem with $CoV=10\%$ (left) and $25\%$
(right) using stochastic collocation (SC), inexact stochastic inverse
subspace iteration (SISI), and inexact line-search Newton method (NI) with the stopping criteria $\Vert
r_{n}\Vert _{2}<10^{-10}$. Here $%
d$ is the polynomial degree and $k$ is the index of basis function in
expansion~(\protect\ref{eq:sol_mat}).}
\label{tab:gpc}
\begin{center}
{\footnotesize \renewcommand{\arraystretch}{1.3} 
\begin{tabular}{|c|c|rrr|rrr|}
\hline
&  & \multicolumn{3}{|c|}{$CoV=10\%$} & \multicolumn{3}{c|}{$CoV=25\%$} \\ 
\hline
$d$ & $k$ & \multicolumn{1}{c}{SC} & \multicolumn{1}{c}{SISI} & \multicolumn{1}{c}{NI} & \multicolumn{1}{|c}{SC} & \multicolumn{1}{c}{SISI} & \multicolumn{1}{c|}{NI} \\ \hline
0 & 1 & 4.9431E+00 & 4.9431E+00 & 4.9431E+00 & 4.9052E+00 & 4.9052E+00 & 
4.9052E+00 \\ \hline
\multirow{3}{*}{1} & 2 & 3.6197E-01 & 3.6197E-01 & 3.6197E-01 & 8.8127E-01 & 
8.8127E-01 & 8.8127E-01 \\ 
& 3 & 1.4477E-13 & -1.6489E-14 & -7.9829E-16 & 2.0162E-13 & -1.5964E-14 & 
-7.3784E-16 \\ 
& 4 & -6.6436E-13 & -1.7135E-14 & -1.3429E-15 & 9.9476E-14 & -1.8588E-14 & 
-1.4099E-15 \\ \hline
\multirow{6}{*}{2} & 5 & 1.8642E-02 & 1.8642E-02 & 1.8642E-02 & 1.1205E-01 & 
1.1201E-01 & 1.1204E-01 \\ 
& 6 & -5.4534E-13 & -9.5178E-17 & -7.4261E-17 & -7.1498E-14 & -2.9421E-15 & 
-1.6150E-16 \\ 
& 7 & -3.0909E-13 & -1.1628E-15 & -9.5249E-17 & -9.4147E-14 & -2.4433E-15 & 
-3.7169E-16 \\ 
& 8 & -1.5442E-03 & -1.5442E-03 & -1.5442E-03 & -9.1479E-03 & -9.1520E-03 & 
-9.1493E-03 \\ 
& 9 & -9.7700E-15 & -1.1200E-15 & 1.3125E-18 & -8.4643E-13 & 7.4442E-16 & 
-1.2278E-17 \\ 
& 10 & -1.5442E-03 & -1.5442E-03 & -1.5442E-03 & -9.1479E-03 & -9.1520E-03 & 
-9.1493E-03 \\ \hline
\end{tabular}
}
\end{center}
\end{table}

A comparison of the inexact SISI\ and the inexact Newton iteration~(NI) is
provided by Figure~\ref{fig:sisi-ni_comparison}, which shows the $2$-norms
of the residual indicator~$\bar{\widetilde{r}}^{s,\left( n\right) }=[%
\widetilde{r}_{1}^{s,(n)T},\dots ,\widetilde{r}_{n_{\xi }}^{s,(n)T}]^{T}$
from~(\ref{eq:res-indicator}) and the part of the residual in the Newton
method given by$~F(\bar{u}^{s,(n)},\bar{\lambda}^{s,(n)})$, cf.~(\ref%
{eq:N-system}). These quantities correspond to the residual of eq.~(\ref%
{eq:SG-eig}), through eq.~(\ref{eq:SG-eig-2}) and equivalent~eq.~(\ref%
{eq:SG-eig-3}). It can be seen that it takes approximately the same number
of steps for the NI to converge and for the SISI residuals to become flat in
case of repeated eigenvalues, but more steps of SISI are needed for simple
eigenvalues. With~respect to the average number of Krylov iterations per a
step of SISI and NI, the computational cost of the two methods is comparable
for simple eigenvalues, but SISI\ is significantly more efficient for
repeated eigenvalues. On the other hand, NI\ outperforms SISI in terms of
accuracy of the solution residual, which is quite natural since NI\ is
formulated as a minimization algorithm unlike SISI.

\begin{figure}[!t]
\vspace{-5mm}
\centering
\subfloat{	\includegraphics[angle=0, scale=0.58]{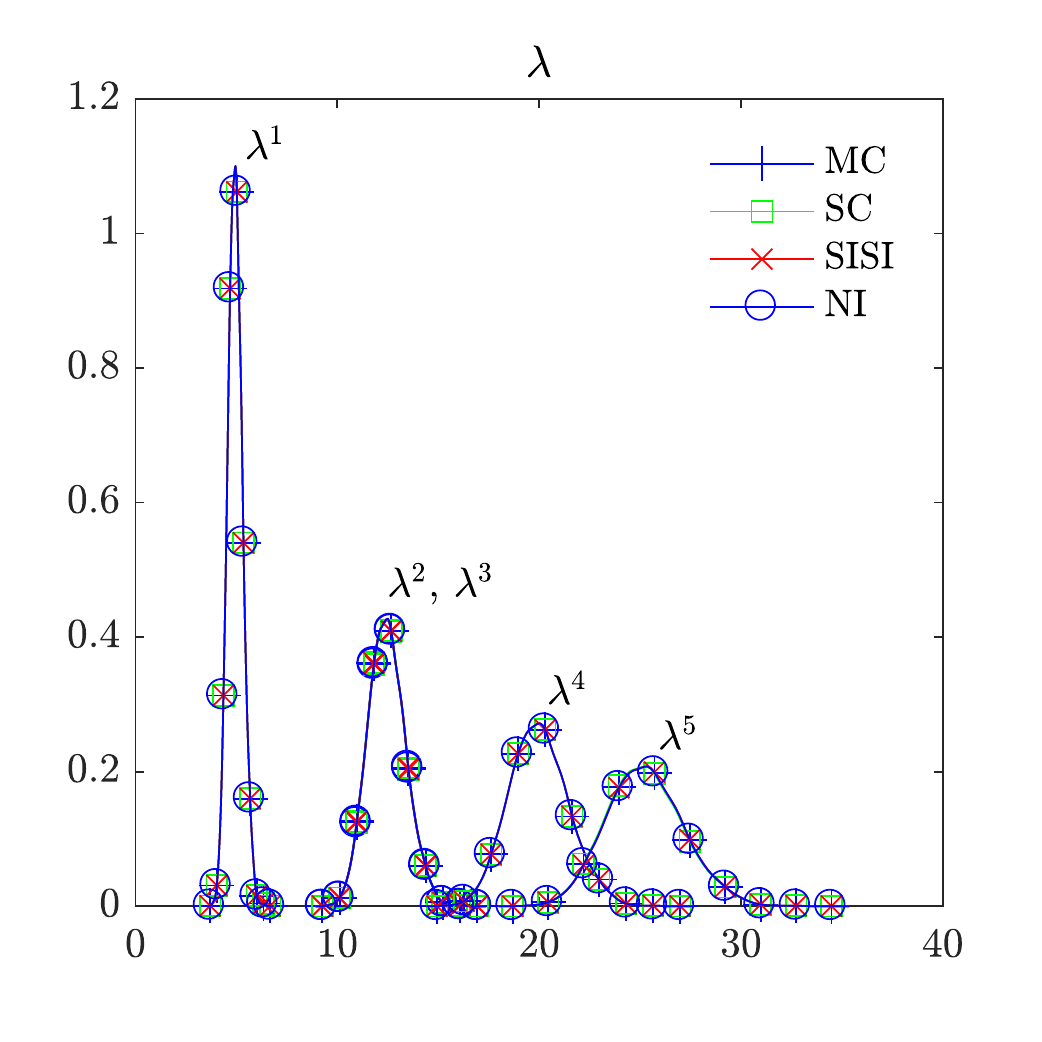}
	\label{fig:eigval_test1}
	} 
\subfloat{	\includegraphics[angle=0, scale=0.58]{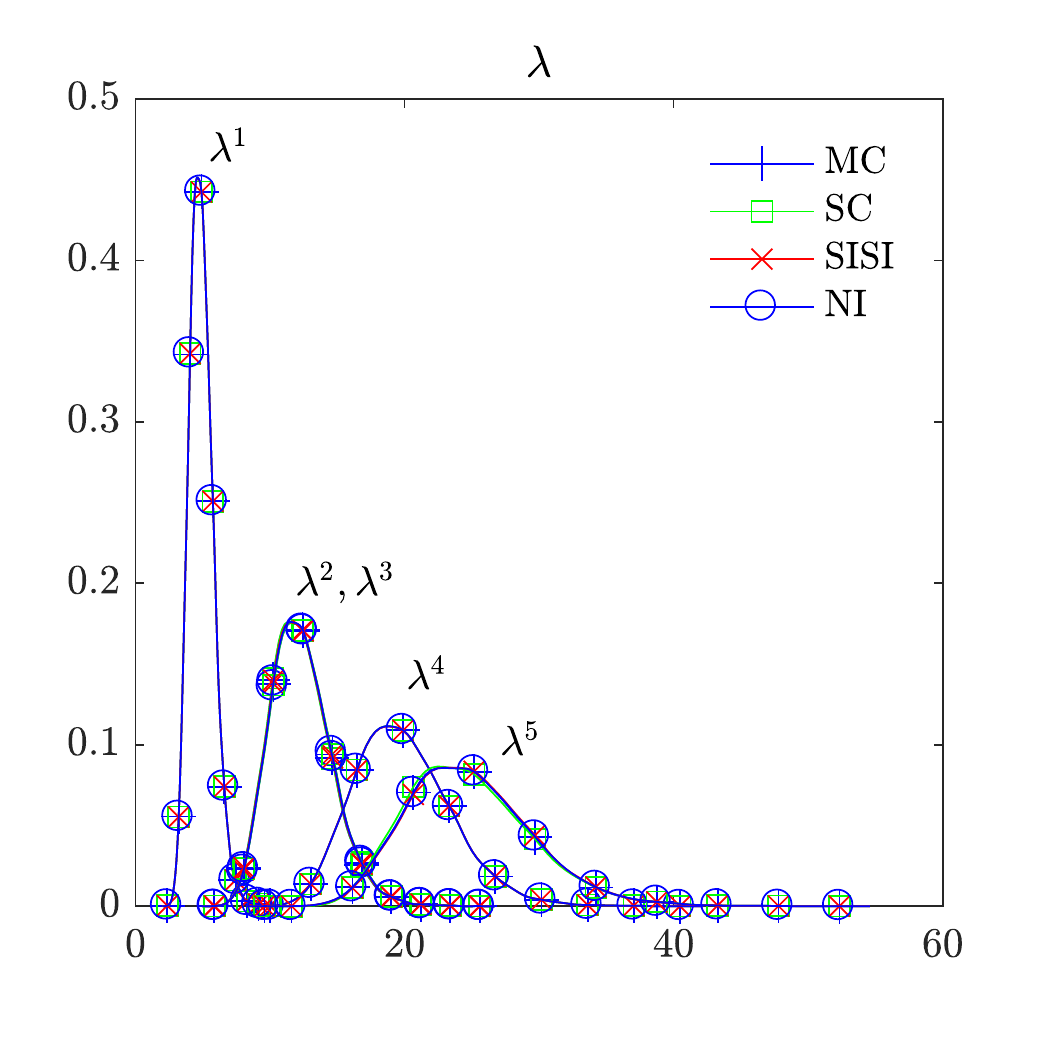}
	\label{fig:eigval_test2}
	}
	     \vspace{-5mm}
\caption{Pdf estimates of the five smallest eigenvalues with $CoV=10\%$
(left) and $25\%$ (right).}
\label{fig:eigvals}
\end{figure}

We also compare the gPC coefficients of eigenvalue expansions computed using
the three different methods: the stochastic collocation method, the inexact
SISI method, and the inexact line-search Newton method. In Table~\ref%
{tab:gpc}, we tabulate the first ten coefficients of the gPC expansion of
the smallest eigenvalues computed using the three methods. A good agreement
of coefficients can be seen, in particular for coefficients with values much
larger than zero, specifically with indices $k=1,2,5,8$ and$~10$. 
Figure~\ref{fig:eigvals} plots the probability density function
(pdf) estimates of the five smallest eigenvalues obtained directly by Monte
Carlo and the three methods, for which the estimates were obtained using 
\textsc{Matlab} function \texttt{ksdensity} used for sampled gPC expansions.
It can be seen that the pdf estimates overlap in all cases. 

\begin{table}[!b]
\caption{Average number of PCG iterations for computing the five smallest
eigenvalues and corresponding eigenvectors of the diffusion problem 
with $CoV =10\%$ (left) and $25\%$ (right) using exact and inexact stochastic inverse
subspace iteration (Algorithm~\protect\ref{alg:isisi}).}
\label{tab:pcg_iters_exact_inexact}
\begin{center}
{\footnotesize \renewcommand{\arraystretch}{1.3} 
\begin{tabular}{|c|c|rrrrr|rrrrr|}
\hline
& & \multicolumn{5}{c|}{$CoV=10\%$} & \multicolumn{5}{c|}{$CoV=25\%$} \\ 
\hline
&  & \multicolumn{1}{c}{1st} & \multicolumn{1}{c}{2nd} & \multicolumn{1}{c}{3rd} & \multicolumn{1}{c}{4th} & \multicolumn{1}{c}{5th} & \multicolumn{1}{|c}{1st} & \multicolumn{1}{c}{2nd} & \multicolumn{1}{c}{3rd} & \multicolumn{1}{c}{4th} & \multicolumn{1}{c|}{5th}
\\ \hline
\multirow{2}{*}{MB} & {\scriptsize Inexact} & 6.45 & 3.90 & 3.90 & 4.60 & 3.75 & 8.60 & 5.55 & 5.55 & 6.05 & 4.75 \\ 
 \cline{2-12}
 & {\scriptsize Exact} &11.00 & 10.95 & 10.95 & 10.85 & 10.00 & 17.00 & 16.90 & 16.90 & 16.90 & 16.90\\
\hline
hGS &{\scriptsize Inexact} & 2.35 & 1.70 & 1.70 & 1.65 & 1.00 & 2.60 & 1.90 & 1.90 & 1.85 & 1.75 \\ 
 \cline{2-12}
$\!\!(p_t\!\!=\!\!2)\!\!$& {\scriptsize Exact} & 3.00 & 3.00 & 3.00 & 3.00 & 3.00 & 5.00 & 5.00 & 5.00 & 4.00 & 4.00\\
\hline
\end{tabular}
}
\end{center}
\end{table}
\begin{table}[!b]
\caption{The number of GMRES iterations for computing the first and the fourth smallest
eigenvalues and corresponding eigenvectors of the diffusion problem with $%
CoV = 10\%$ (left) and $25\%$ (right) using exact and inexact line-search Newton
method (Algorithm~\protect\ref{alg:line_search}) with preconditioners cMB
(top) and chGS($p_t=2$) (bottom), and with the stopping criteria $\Vert
r_{n}\Vert _{2}<10^{-10}$.}
\label{tab:gmres_exact_inexact}
\begin{center}
{\footnotesize \renewcommand{\arraystretch}{1.3} 
\begin{tabular}{|c|cc|cc||cc|cc|}
\hline
& \multicolumn{4}{|c||}{$CoV=10\%$} & \multicolumn{4}{c|}{$CoV=25\%$} \\ 
\hline
& \multicolumn{2}{|c|}{Inexact} & \multicolumn{2}{c||}{Exact} & \multicolumn{2}{c|}{Inexact} & \multicolumn{2}{c|}{Exact}  \\ 
\hline
& 1st & 4th & 1st & 4th & 1st & 4th & 1st & 4th  \\ \hline
Nonlinear step & \multicolumn{8}{|c|}{cMB{\tiny {\ (updated)}}} \\ \hline
1 & 2 & 1   & 13 & 16 & 2   & 1   & 22 & 39  \\ 
2 & 4 & 3   & 13 & 15 & 5   & 3   & 22 & 27 \\ 
3 & 7 & 6   & 14 & 16 & 9   & 6   & 22 & 27  \\ 
4 &    & 11 &      & 16 & 15 & 12 & 22 & 27  \\ 
5 &    &      &     &       &      & 21 &      & 27\\ 
\hline
Nonlinear step & \multicolumn{8}{|c|}{chGS($p_t=2$)} \\ \hline
1 & 1 & 1 & 5 & 7 & 1 & 1 & 7 & 16  \\ 
2 & 2 & 1 & 5 & 7 & 2 & 1 & 8 & 15 \\ 
3 & 3 & 2 & 6 & 7 & 4 & 2 & 8 & 11 \\ 
4 &    & 4 &    & 7 & 6 & 3 & 8 & 11 \\ 
5 &    &   &    &     &    & 5 &    & 10 \\ 
6 &    &   &    &     &    &    &    & 10 \\ 
\hline
\end{tabular}
}
\end{center}
\end{table}

\subparagraph{Inexact vs. exact solves}
We present numerical experiments that show the effectiveness of the inexact solvers by comparing them with the exact solvers, for which we fix the stopping tolerance of the PCG and GMRES methods to $10^{-12}$.
For the inexact methods we use the adaptive stopping tolerance given for SISI by~(\ref{eq:pcg-stop})
and for the NI by~(\ref{eq:gmres-stop}). 
A comparison of the inexact and exact solves in terms of the PCG iteration counts 
for computing the smallest five eigenvalues of the diffusion problem is shown in Table~\ref{tab:pcg_iters_exact_inexact},
and a comparison in terms of the GMRES iterations counts for computing the first and the fourth smallest eigenvalues of the diffusion problem 
is shown in Table~\ref{tab:gmres_exact_inexact}. 
In both cases, for given $CoV$ and the choice of the preconditioner, we observe that the exact methods require more Krylov subspace iterations. 
It can be seen from Table~\ref{tab:gmres_exact_inexact} that virtually the same number of GMRES
iterations is required in each nonlinear step of NI since the stopping tolerance of the exact solves is not adjusted to the nonlinear residual. 

\subparagraph{Effect of increasing the stochastic dimension}
Table~\ref{tab:pcg_iters_mxi} shows the PCG iteration counts required to compute the smallest five eigenvalues of the diffusion problem for varying 
number of random variables $m_{\xi}=\{3,5,7\}$ with $CoV = 10\%$ and $25\%$, 
and Table \ref{tab:lognormal_higher_dims} shows the GMRES iteration counts 
for computing the first and fourth smallest eigenvalues for the same problem and setup. 
While in both cases we see a relatively small increase in iteration counts for larger $CoV$, 
increasing the stochastic dimension by setting larger~$m_{\xi}$ appears to have no effect on the iteration counts.

\begin{table}[!h]
\caption{Average number of PCG iterations for computing the five smallest
eigenvalues and corresponding eigenvectors of the diffusion problem with $%
CoV =10\%$ (left) and $25\%$ (right) for varying $m_{\xi}=\{3,5,7\}$ using inexact stochastic inverse
subspace iteration (Algorithm~\protect\ref{alg:isisi}).}
\label{tab:pcg_iters_mxi}
\begin{center}
{\footnotesize \renewcommand{\arraystretch}{1.3} 
\begin{tabular}{|c|c|ccccc|ccccc|}
\hline
&& \multicolumn{5}{c|}{$CoV=10\%$} & \multicolumn{5}{c|}{$CoV=25\%$} \\ 
\hline
$m_{\xi}$&Preconditioner & 1st & 2nd & 3rd & 4th & 5th & 1st & 2nd & 3rd & 4th & 5th
\\ \hline
\multirow{2}{*}{3}&MB & 6.45 & 3.90 & 3.90 & 4.60 & 3.75 & 8.60 & 5.55 & 5.55 & 6.05 & 4.75 \\ 
&hGS ($p_t=2$) & 2.35 & 1.70 & 1.70 & 1.65 & 1.00 & 2.60 & 1.90 & 1.90 & 1.85
& 1.75 \\ 
 \hline
 \multirow{2}{*}{5}&MB & 6.50 & 3.90 & 3.90 & 4.50 & 3.85 & 8.00 & 4.85 & 4.85 & 6.50 & 4.70 \\ 
&hGS ($p_t=2$) & 2.35 & 1.00 & 1.00 & 1.70 & 1.00 & 2.60 & 1.95 & 1.95 & 1.90 & 1.85\\
 \hline
 \multirow{2}{*}{7}&MB & 6.40 & 3.95 & 3.95 & 4.55 & 3.85 & 8.00 & 4.85 & 4.85 & 6.50 & 4.70 \\ 
&hGS ($p_t=2$) &2.35 & 1.00 & 1.00 & 1.70 & 1.00 & 2.60 & 1.95 & 1.95 & 1.90 & 1.85\\
 \hline
\end{tabular}
}
\end{center}
\end{table}

\begin{table}[!t]
\caption{The number of GMRES iterations for computing the first and the fourth smallest
eigenvalues and corresponding eigenvectors of the diffusion problem with $%
CoV = 10\%$ (left) and $25\%$ (right) for varying $m_{\xi}$ using inexact line-search Newton
method (Algorithm~\protect\ref{alg:line_search}) with preconditioners cMB
(top) and chGS($p_t=2$) (bottom), and with the stopping criteria $\Vert
r_{n}\Vert _{2}<10^{-10}$.}
\label{tab:lognormal_higher_dims}
\begin{center}
{\footnotesize \renewcommand{\arraystretch}{1.3} 
\begin{tabular}{|c|cc|cc|cc||cc|cc|cc|}
\hline
& \multicolumn{6}{c||}{$CoV=10\%$} & \multicolumn{6}{c|}{$CoV=25\%$} \\ 
\hline
$m_{\xi}$& \multicolumn{2}{c|}{3} & \multicolumn{2}{|c|}{5} & \multicolumn{2}{c||}{7} & \multicolumn{2}{c|}{3} & \multicolumn{2}{c|}{5} & \multicolumn{2}{|c|}{7}  \\ 
\hline
& 1st & 4th & 1st & 4th & 1st & 4th & 1st & 4th  & 1st & 4th & 1st & 4th  \\ \hline
Nonlinear step & \multicolumn{12}{c|}{cMB{\tiny {\ (updated)}}} \\ \hline
1 & 2 & 1   & 2 & 1 & 2   & 1   & 2    & 1  & 2   & 1  & 2  & 1\\ 
2 & 4 & 3   & 4 & 3 & 4   & 3   & 5    & 3  & 5   & 3  & 4  & 3\\ 
3 & 7 & 6   & 7 & 6 & 7   & 6   & 9    & 6  & 9   & 6  & 8  & 6\\ 
4 &    & 11 &    & 11 &    & 11  & 15  & 12& 15 & 12&16 & 12\\  
5 &    &      &     &    &     &      &        & 21&     &  21&    & 21\\
\hline
Nonlinear step & \multicolumn{12}{c|}{chGS($p_t=2$)} \\ \hline
1 & 1 & 1 & 1 & 1 & 1 & 1 & 1 & 1 & 1 & 1 & 1 & 1 \\ 
2 & 2 & 1 & 2 & 1 & 2 & 1 & 2 & 1 & 2 & 1 & 2 & 1 \\ 
3 & 3 & 2 & 3 & 2 & 3 & 2 & 4 & 2 & 3 & 2 & 4 & 2\\ 
4 &    & 4 &    & 4 &    & 4 & 6 & 3 & 5 & 3 & 5 & 3 \\ 
5 &    &   &    &     &    &    &    & 5 &    & 5 &    &  5\\ 
\hline
\end{tabular}
}
\end{center}
\end{table}

\subsection{Stiffness of Mindlin plate with uniformly distributed Young's modulus}
As the second example, we study eigenvalues of the stiffness of Mindlin plate
with Young's modulus given by the stochastic expansion  
\begin{equation}\label{eq:mindlin_exp}
E(x,\xi) = E_1 + \sum_{\ell=2}^{m_{\xi}+1} E_\ell \xi_{\ell-1},
\end{equation}
where $E_{\ell+1}=\sqrt{\lambda_\ell} v_\ell(x)$ with $\{ \left( \lambda_\ell,v_\ell \right) \}_{\ell=1}^{m_{\xi}}$ are the eigenpairs 
of the eigenvalue problem associated with the covariance kernel
\begin{equation}
C\left(  X_{1},X_{2}\right)  = \frac{1}{3} \sigma_{u}^{2}\exp\left(  -\frac{\left\vert
x_{2}-x_{1}\right\vert }{L_{x}}-\frac{\left\vert y_{2}-y_{1}\right\vert
}{L_{y}}\right)  , \label{eq:covariance-2}%
\end{equation}
where $L_{x}$, $L_{y}$ are as in~(\ref{eq:covariance}),
and $\sigma_{u}$ is the standard deviation of the random field,
the random variables $\xi_\ell$ are uniformly distributed over the interval $(-1,1)$, $E_1=10920$, 
and other parameters are set as in~\cite{Sousedik-2016-ISI}. 
The plate is discretized using $10 \times 10$ bilinear (Q4) finite elements with $243$ physical degrees of freedom. 
We note that we consider only the stiffness matrix in the problem setup, and the mass matrix is taken as identity.
For the uniform random variables, the set~$\{\psi_{k}\}_{k=1}^{n_{\xi}}$ is given by Legendre polynomials and   
Smolyak sparse grid with Gauss-Legendre quadrature is considered for the quadrature rule. % and grid level~$4$,

Table~\ref{tab:pcg_iters_Mindlin} shows the average numbers of PCG iterations required 
to solve linear system~(\ref{eq:alg-SISI-solve})
with zero initial guess and the adaptive stopping criteria~(\ref{eq:pcg-stop}). 
As we observed in the results of the diffusion problem in Table~\ref{tab:pcg_iters}, 
PCG with the hGS preconditioning requires less than the half of the iteration counts with the MB preconditioner.
Table~\ref{tab:mindlin_NI} shows the average numbers of GMRES iterations required to solve 
the linear systems at Line~\ref{ln:lin_sys} in Algorithm~\ref{alg:line_search} with zero initial guess 
and the adaptive stopping criteria~(\ref{eq:gmres-stop}).  
As in the results of the diffusion problem in Table~\ref{tab:nl_iters}, we again observe that  the updated versions 
of the preconditioners yield lower iteration counts compared to their fixed variants and the lowest counts are achieved 
with the chGS preconditioner. Increasing both $CoV$ and stochastic dimension $m_{\xi}$ leads to only a mild increase 
in iteration counts. Finally, Table~\ref{tab:gpc_Mindlin} shows the first $10$ coefficients of the gPC expansion 
of the smallest eigenvalue of the Mindlin plate. 
As for the solution coefficients of the diffusion problem shown in Table~\ref{tab:gpc}, 
a good agreement of coefficients can be seen also here. 

\begin{table}[!t]
\caption{Average number of PCG iterations for computing the five smallest
eigenvalues and corresponding eigenvectors of the Mindlin plate problem with $%
CoV =10\%$ (left) and $25\%$ (right) using inexact stochastic inverse
subspace iteration (Algorithm~\protect\ref{alg:isisi}).}
\label{tab:pcg_iters_Mindlin}
\begin{center}
{\footnotesize \renewcommand{\arraystretch}{1.3} 
\begin{tabular}{|c|ccccc|ccccc|}
\hline
& \multicolumn{5}{c|}{$CoV=10\%$} & \multicolumn{5}{c|}{$CoV=25\%$} \\ 
\hline
Preconditioner & 1st & 2nd & 3rd & 4th & 5th & 1st & 2nd & 3rd & 4th & 5th
\\ \hline
MB & 6.20	 & 4.65 & 4.65 & 4.70  & 4.20 & 8.15 & 6.55 & 6.55 & 6.75 & 6.05\\ 
hGS ($p_t=1$) & 2.45 & 1.95 & 1.95 & 1.95 & 1.95 & 3.40 & 2.75 & 2.75 & 2.65 & 2.60\\ 
hGS ($p_t=2$) & 2.45 & 1.95 & 1.95 & 1.95 & 1.95 & 3.40 & 2.75 & 2.75 & 2.65 & 2.60\\ 
hGS (no trunc.) & 2.45 & 1.95 & 1.95 & 1.95 & 1.95 & 3.40 & 2.75 & 2.75 & 2.65 & 2.60\\ 
\hline
\end{tabular}
}
\end{center}
\end{table}

\begin{table}[!t]
\caption{The average number of GMRES iterations for computing the first and the fourth smallest
eigenvalues and corresponding eigenvectors of the Mindlin plate problem with $%
CoV = 10\%$ and  $25\%$  for varying $m_{\xi}$ (the number of random variables) using inexact line-search Newton
method (Algorithm~\protect\ref{alg:line_search}) with preconditioners cMB
(top) and chGS($p_t=2$) (bottom), and with the stopping criteria $\Vert
r_{n}\Vert _{2}<10^{-10}$.}
\label{tab:mindlin_NI}
\begin{center}
{\footnotesize \renewcommand{\arraystretch}{1.3} 
\begin{tabular}{|c|c|rr|rr|rr|rr|}
\hline
\multirow{2}{*}{}& $m_{\xi}$ &  \multicolumn{2}{c|}{$3$} & \multicolumn{2}{|c|}{$5$}& \multicolumn{2}{c|}{$7$}& \multicolumn{2}{|c|}{$9$}\\
\cline{2-10}
& & 1st & 4th & 1st & 4th & 1st & 4th& 1st & 4th  \\
\hline
\multirow{8}{*}{$CoV=10\%$}&NMB{\tiny {\ (fixed)}}&  14.25 & 26.50 & 15.25 & 30.75 & 15.25 & 33.25 & 15.25 & 34.00\\
\cline{2-10}
&NMB{\tiny {\ (updated)}} & 12.00 & 12.00 & 15.00 & 13.75 & 15.00	& 14.00 & 15.00 & 14.25\\
\cline{2-10}
& cMB{\tiny {\ (fixed)}} & 10.25 & 10.25 & 10.75 & 11.25 & 11.00 & 11.50 & 11.00 & 11.75\\
\cline{2-10}
& cMB{\tiny {\ (updated)}} & 6.00 & 5.25 & 6.25 & 5.75 & 6.25 & 6.00 & 6.75 & 6.00\\
%\cline{2-10}
%& chGS($p_t=0$) & 6.00 & 5.25 & 6.25 & 5.75 & 6.25 & 6.00 & 6.75 & 6.00\\
\cline{2-10}
& chGS($p_t=1$) & 3.00 & 2.75 & 3.00 & 3.00 & 3.00 & 3.00 & 3.00 & 3.00\\
\cline{2-10}
& chGS($p_t=2$) & 3.00 & 2.75 & 3.00 & 2.75 & 3.00 & 3.00 & 3.00 & 3.00\\
\cline{2-10}
& chGS(full) & 3.00 & 2.75 & 3.00 & 2.75 & 3.00 & 3.00 & 3.00 & 3.00\\
 \hline
\multirow{8}{*}{$CoV=25\%$} &NMB{\tiny {\ (fixed)}} & 13.25 & 32.40 & 14.50 & 42.80 & 14.75 & 61.20 & 20.00 & 63.60\\
\cline{2-10}
& NMB{\tiny {\ (updated)}} & 14.75 & 16.60 & 19.75 & 29.17 & 26.40	& 40.00 & 27.60 & 42.67\\
\cline{2-10}
& cMB{\tiny {\ (fixed)}} & 11.25 & 18.17 & 12.50 & 22.33 & 12.50 & 28.83 & 17.40 & 29.50\\
\cline{2-10}
& cMB{\tiny {\ (updated)}} & 6.50 & 10.83 & 7.25 & 12.67 & 10.20 & 16.33 & 10.20 & 17.00\\
%\cline{2-10}
%& chGS($p_t=0$)& 6.50 & 10.83 & 7.25 & 12.67 & 10.20 & 16.33 & 10.20 & 17.00\\
\cline{2-10}
& chGS($p_t=1$)& 3.25 & 4.83 & 3.25 & 5.33 & 3.25 & 7.17 & 4.60 & 7.67\\
\cline{2-10}
& chGS($p_t=2$) & 3.25 & 4.83 & 3.25 & 5.33 & 3.25 & 7.17 & 4.40 & 7.50\\
\cline{2-10}
& chGS(full) & 3.25 & 4.83 & 3.25 & 5.50 & 3.25 & 6.83 & 4.40 & 7.33\\
 \hline
\end{tabular}
}
\end{center}
\end{table}

\begin{table}[!t]
\caption{The first $10$ coefficients of the gPC expansion of the smallest
eigenvalue of the Mindlin plate problem with $CoV=10\%$ (left) and $25\%$
(right) using stochastic collocation (SC), inexact stochastic inverse
subspace iteration (SISI), and inexact line-search Newton method (NI) with the stopping criteria $\Vert
r_{n}\Vert _{2}<10^{-10}$. Here $%
d$ is the polynomial degree and $k$ is the index of basis function in
expansion~(\protect\ref{eq:sol_mat}).}
\label{tab:gpc_Mindlin}
\begin{center}
{\footnotesize \renewcommand{\arraystretch}{1.3} 
\begin{tabular}{|c|c|rrr|rrr|}
\hline
&  & \multicolumn{3}{c|}{$CoV=10\%$} & \multicolumn{3}{c|}{$CoV=25\%$} \\ 
\hline
$d$ & $k$ & \multicolumn{1}{c}{SC} & \multicolumn{1}{c}{SISI} & \multicolumn{1}{c}{NI} &\multicolumn{1}{|c}{SC} & \multicolumn{1}{c}{SISI} & \multicolumn{1}{c|}{NI} \\ \hline
0 & 1 & 4.6271E-01 & 4.6271E-01 & 4.6271E-01 & 4.5784E-01 & 4.5784E-01 & 4.5784E-01 \\ \hline
\multirow{3}{*}{1} & 2 & -2.2476E-02 & -2.2476E-02 & -2.2476E-02 & -5.6737E-02 & -5.6734E-02 & -5.6735E-02 \\ 
& 3 & 6.6391E-14 & -3.5389E-16 & -8.0416E-18 & -1.7453E-13 & -6.5624E-16 & -1.1174E-17\\ 
& 4 & 3.2080E-13 & -4.2037E-16 & 1.4672E-17	 & 6.0396E-14 & -4.8016E-16 & 2.5675E-17\\ 
\hline
\multirow{6}{*}{2} & 5 & -3.1659E-05 & -3.1607E-05 & -3.1634E-05 & -2.5953E-04 & -2.4582E-04 & -2.5268E-04\\ 
& 6 & -7.8920E-14 & 1.7146E-16 & -1.0762E-18 & -2.2204E-16 & 9.0132E-16 & 4.9237E-18\\ 
& 7 & 3.1186E-13 & 3.8511E-16 & -4.5709E-19 & -6.1270E-15 & 9.5916E-16 & 8.0412E-18\\ 
& 8 & -3.8995E-04 & -3.8995E-04 & -3.8995E-04 & -2.5032E-03 & -2.5021E-03 & -2.5030E-03\\ 
& 9 & -2.8144E-14 & -9.5150E-17 & -5.8430E-19 & 1.1297E-13 & -9.2077E-17 & -1.5950E-18\\ 
& 10 & -3.8995E-04 & -3.8995E-04 & -3.8995E-04 & -2.5032E-03 & -2.5021E-03 & -2.5030E-03\\
\hline
\end{tabular}
}
\end{center}
\end{table}

\section{Conclusion}

\label{sec:conclusion}We studied inexact methods for symmetric eigenvalue
problems in the context of spectral stochastic finite element
discretizations. The performance was compared using eigenvalue problems
given by the stochastic diffusion equation with lognormally distributed
diffusion coefficient and by the stiffness of Mindlin plate with Young's modulus depending 
on uniformly distributed random variables. 
Both problems were given in a $2$-dimensional physical domain. The methods were
formulated on the basis of the stochastic inverse subspace iteration (SISI)
and the line-search Newton method (NI). In both formulations we obtained
symmetric stochastic Galerkin matrices. In the first case the matrices were
also positive definite, so the associated linear systems were solved using
preconditioned conjugate gradient (PCG) method. For the PCG we used
mean-based and hierarchical Gauss-Seidel preconditioners. The second
preconditioner slightly decreased the overall iteration count, but in all
cases only a handful of iterations were required for convergence per one
step of SISI. The iteration count for PCG also did not appear to be
sensitive to algebraic multiplicity of eigenvalues, but in terms of SISI we
observed somewhat slower convergence for simple eigenvalues (i.e., those
with algebraic multiplicity one). For the second method based on Newton
iteration, we proposed several novel preconditioners adapted to the
structure of the Jacobian matrices obtained from the stochastic Galerkin
discretization. The linear systems were solved using the GMRES (and in a few
cases also MINRES) method with various preconditioners. 
We analytically show that chGS with a truncated matrix-vector product is 
the most efficient one for high-dimensional problems. 
The overall iteration count of GMRES was higher compared to PCG, 
in particular for eigenvalues with algebraic multiplicity
larger than one. On the other hand, only a handful of iterations were
required with the constraint hierarchical Gauss-Seidel preconditioner for
simple eigenvalues. In terms of the iteration count of the SISI and NI, we
observed that the two methods are comparable for simple eigenvalues, but
SISI appeared more efficient for repeated eigenvalues. Increasing either the value
of$~CoV$ or the stochastic dimension lead to only a slight increase of the number of iterations, in
particular when the constraint hierarchical preconditioners were used.
Comparing the accuracy in terms of the solution residual,
NI\ naturally outperformed SISI.\ Nevertheless both methods identified the
coefficients of polynomial chaos expansion of the smallest eigenvalue in a
close agreement\ and matched well those computed by the stochastic
collocation. The probability density estimates of all eigenvalues matched,
also with the direct Monte Carlo simulation.

From a user's perspective, the SISI is straightforward to use and in
combination with the stochastic modified Gram-Schmidt process allows to
compute coefficients of polynomial chaos expansions of several eigenvalues
and eigenvectors, while the NI\ requires some setup of parameters for the
line search and backtracking. On the other hand, NI may be more suitable
when interior eigenvalues are sought, since the SISI assumes that all
smaller eigenvalues were deflated from the mean matrix.

\paragraph{Acknowledgement} We would like to thank Prof. Howard C. Elman for 
sharing his pearls of wisdom with us and many fruitful discussions. 
We would also like to thank the anonymous referees for insightful comments. 
This paper describes objective technical results and analysis. Any subjective
views or opinions that might be expressed in the paper do not necessarily
represent the views of the U.S. Department of Energy or the United States
Government. Sandia National Laboratories is a multimission laboratory managed
and operated by National Technology and Engineering Solutions of Sandia, LLC.,
a wholly owned subsidiary of Honeywell International, Inc., for the U.S.
Department of Energy's National Nuclear Security Administration under contract
DE-NA-0003525.

\appendix

\section{Inexact Newton iteration}
\label{sec:ini}The inexact nonlinear iteration is based on the
Newton--Krylov method, in which each step entails solving the linear system~(%
\ref{eq:Newton-mod}) by a Krylov subspace method followed by an update~(\ref%
{eq:Newton-update}). %For simplicity, we drop the superscript~$^{s,(n)}$ 
%in the rest of this section. 
But first, let us describe the evaluation of~$F(\bar{u}^{s,(n)},\bar{\lambda}%
^{s,(n)})$ and~$G(\bar{u}^{s,(n)})$. The vector$~F(\bar{u}^{s,(n)},\bar{%
\lambda}^{s,(n)})$, defined by~(\ref{eq:F_mat}), consists of two terms: the
first term is evaluated as 
\begin{equation*}
\mathbb{E}[\Psi \Psi ^{T}\otimes A]\bar{u}^{s,(n)}=\sum_{\ell =1}^{n_{a}}({H}%
_{\ell }\otimes A_{\ell })\bar{u}^{s,(n)}=\text{vec}\left( \sum_{\ell
=1}^{n_{a}}A_{\ell }\bar{U}^{s,(n)}H_{\ell }^{T}\right) ,
\end{equation*}%
which is the same as~(\ref{eq:mat-vec}), and the second term is evaluated as 
\begin{equation*}
\mathbb{E}[((\bar{\lambda}^{s,(n)})^{T}\Psi )\Psi \Psi ^{T}\otimes
I_{n_{x}}]\bar{u}^{s,(n)}=\sum_{i=1}^{n_{\xi }}({\lambda}%
_{i}^{s,(n)}H_{i}\otimes I_{n_{x}})\bar{u}^{s,(n)}=\text{vec}\left(
\sum_{i=1}^{n_{\xi }}{\lambda}_{i}^{s,(n)}\bar{U}^{s,(n)}H_{i}^{T}%
\right) .
\end{equation*}%
The vector$~G(\bar{u}^{s,(n)})$, defined by~(\ref{eq:G_mat}), is evaluated
as 
\begin{equation*}
G(\bar{u}^{s,(n)})=\mathbb{E}\left[ \Psi \otimes \left( (\bar{u}%
^{s,(n)}{}^{T}(\Psi \Psi ^{T}\otimes I_{n_{x}})\bar{u}^{s,(n)})-1\right) %
\right] ,
\end{equation*}%
where the $i$th row of$~G(\bar{u}^{s,(n)})$ is 
\begin{align*}
\left[ G(\bar{u}^{s,(n)})\right] _{i}& =\mathbb{E}[\psi _{i}(\bar{u}%
^{s,(n)}{}^{T}(\Psi \Psi ^{T}\otimes I_{n_{x}})\bar{u}^{s,(n)})-\psi _{i}],
\\
& =\bar{u}^{s,(n)}{}^{T}\mathbb{E}[\psi _{i}\Psi \Psi ^{T}\otimes I_{n_{x}}]%
\bar{u}^{s,(n)}-\delta _{1i},
\end{align*}%
and the first term above is evaluated as 
\begin{equation*}
\bar{u}^{s,(n)}{}^{T}\mathbb{E}[\psi _{i}\Psi \Psi ^{T}\otimes I_{n_{x}}]%
\bar{u}^{s,(n)}=\bar{u}^{s,(n)}{}^{T}(H_{i}\otimes I_{n_{x}})\bar{u}^{s,(n)},
\end{equation*}%
or, denoting the trace operator by$~\text{tr}$, this term can be also
evaluated as 
\begin{equation*}
\bar{u}^{s,(n)}{}^{T}\mathbb{E}[\psi _{i}\Psi \Psi ^{T}\otimes I_{n_{x}}]%
\bar{u}^{s,(n)}=\text{tr}(\bar{U}^{s,(n)}H_{i}\bar{U}^{s,(n)}{}^{T})=\text{tr%
}(\bar{U}^{s,(n)}{}^{T}\bar{U}^{s,(n)}H_{i}).
\end{equation*}
\begin{remark}
\label{rem:explicit_jac}For completeness, let us describe a possible setup
of the Jacobian matrices in~(\ref{eq:Newton}) or~(\ref{eq:Newton-mod}).
Block~(\ref{eq:jac_Fu}) can be set up~as 
\begin{equation}
\mathbb{E}[\Psi \Psi ^{T}\otimes A]-\mathbb{E}[((\bar{\lambda}%
^{s,(n)})^{T}\Psi )\Psi \Psi ^{T}\otimes I_{n_{x}}]=\sum_{i=1}^{n_{a}}{H}%
_{i}\otimes A_{i}-\sum_{i=1}^{n_{\xi }}({\lambda}_{i}^{s,(n)}H_{i}%
\otimes I_{n_{x}}).  \label{eq:jac_Fu-impl}
\end{equation}%
Block~(\ref{eq:jac_Fl}) can be set up as 
\begin{equation}
\mathbb{E}[\Psi ^{T}\otimes (\Psi \Psi ^{T}\otimes I_{n_{x}})\bar{u}%
^{s,(n)}]=\mathbb{E}[(\psi _{1}\Psi \Psi ^{T}\otimes I_{n_{x}})\bar{u}%
^{s,(n)},\ldots ,(\psi _{n_{\xi }}\Psi \Psi ^{T}\otimes I_{n_{x}})\bar{u}%
^{s,(n)}],  \label{eq:jac_Fl-impl}
\end{equation}%
and the$~i$th column of this block is 
\begin{equation}
\mathbb{E}[(\psi _{i}\Psi \Psi ^{T}\otimes I_{n_{x}})\bar{u}%
^{s,(n)}]=(H_{i}\otimes I_{n_{x}})\bar{u}^{s,(n)}=\text{vec}\left( \bar{U}%
^{s,(n)}H_{i}^{T}\right) .  \label{eq:jac_Fl-impl-2}
\end{equation}%
Finally, block~(\ref{eq:jac_Gu}) is the transpose of~(\ref{eq:jac_Fl})
scaled by a factor of $-2$, cf.~(\ref{eq:Newton-mod}).
\end{remark}

In implementation, the explicit setup described in Remark~\ref%
{rem:explicit_jac} is avoided because Krylov subspace methods require only
matrix-vector products. Let us write a product with Jacobian matrix from~(%
\ref{eq:Newton-mod}) at step$~n$ of the nonlinear iteration as 
\begin{equation}
J(\bar{u}^{s,(n)},\bar{\lambda}^{s,(n)})%
\begin{bmatrix}
\delta \bar{u} \\ 
\delta \bar{\lambda}%
\end{bmatrix}%
,\qquad \text{where\ }\:J(\bar{u}^{s,(n)},\bar{\lambda}^{s,(n)})=%
\begin{bmatrix}
A & B^{T} \\ 
B & 0%
\end{bmatrix}%
,  \label{eq:jac-scheme}
\end{equation}%
with~$A$ and~$B^{T}$ denoting the matrices in~(\ref{eq:jac_Fu}) and~(\ref%
{eq:jac_Fl}), respectively. Then, 
%the submatrix-vector products, cf.~(\ref{eq:mat-vec}), are 
\begin{align}
A\,\delta \bar{u}&\!=\!\left( \sum_{\ell =1}^{n_{a}}H_{\ell }\!\otimes\! A_{\ell
}\!-\!\!\sum_{i=1}^{n_{\xi }}H_{i}\!\otimes\! {\lambda}_{i}^{s,(n)}I_{n_{x}}\!
\!\right)\! \delta \bar{u}\!=\!\text{vec}\!\left( \sum_{\ell =1}^{n_{a}}A_{\ell
}\delta \bar{U}H_{\ell }^{T}\!-\!\!\sum_{i=1}^{n_{\xi }}{\lambda}%
_{i}^{s,(n)}\delta \bar{U}H_{i}^{T}\!\right)\!,  \label{eq:jac-A-mult} \\
B^{T}\,\delta \bar{\lambda}&\!=\!-\sum_{i=1}^{n_{\xi }}\delta \lambda
 _{i}\mathbb{E}[\Psi ^{T}\!\otimes\! (\Psi \Psi ^{T}\!\otimes\! I_{n_{x}})]%
\bar{u}^{s,(n)}\!=\!-\!\text{vec}\left( \sum_{i=1}^{n_{\xi }}\delta \lambda
 _{i}\bar{U}^{s,(n)}H_{i}^{T}\right) ,  \label{eq:jac-B^T-mult}
\end{align}
and
\begin{align}
\label{eq:jac-B-mult} 
B\,\delta \bar{u}& =-\mathbb{E}[\Psi \otimes (\bar{u}^{s,(n)}{}^{T}(\Psi
\Psi ^{T}\otimes I_{n_{x}}))]\delta \bar{u}  
 =-%
\begin{bmatrix}
\bar{u}^{s,(n)}{}^{T}(H_{1}\otimes I_{n_{x}})\delta \bar{u} \\ 
\vdots \\ 
\bar{u}^{s,(n)}{}^{T}(H_{n_{\xi }}\otimes I_{n_{x}})\delta \bar{u}%
\end{bmatrix}%
,
\end{align}%
where the $i$th row can be equivalently evaluated as $\text{tr}(H_{i}^{T}%
\bar{U}^{s,(n)}{}^{T}\delta \bar{U})$.

\section{Matrix-vector product in the chGS preconditioner}
\label{sec:mvp}The matrix-vector product with subblocks of the stochastic
Jacobian matrices are performed as in~(\ref{eq:jac-scheme})--(\ref{eq:jac-B-mult}). For example, the matrix-vector product with a subblock of
the $A$-part of the Jacobian matrix, cf.~(\ref{eq:jac-A-mult}), can be written as 
\begin{eqnarray}
\sum_{t\in \mathcal{I}_{t}}(\left[ h_{t,(\ell )(k)}\right] \otimes
A_{t})v_{(k)}^{s} &=&\text{vec}\left( \sum_{t\in \mathcal{I}%
_{t}}A_{t}V_{(k)}^{s}\left[ h_{t,(k)(\ell )}\right] \right) ,
\label{eq:matvec-block-N1} \\
\sum_{t\in \mathcal{I}_{t}}(\left[ h_{t,(\ell )(k)}\right] \otimes \lambda
_{t}^{s,(n)}I_{n_{x}})v_{(k)}^{s} &=&\text{vec}\left( \sum_{t\in \mathcal{I}%
_{t}}\lambda _{t}^{s,(n)}V_{(k)}^{s}\left[ h_{t,(k)(\ell )}\right]
^{T}\right) ,  \label{eq:matvec-block-N2}
\end{eqnarray}%
where$~V_{(k)}^{s}$ is a subset of the columns of$~V^{s}$ specified by the
index set$~(k)$. We note that the matrix-vector products in~(\ref%
{eq:matvec-block-N2}) depend on the eigenvalue approximation at step~$n$ of
Newton iteration. The truncation of the matrix-vector products, indicated by
summing up over index set$~\mathcal{I}_{t}$ is performed using the same
strategy as in Algorithm~\ref{alg:hGS}. 

\section{Computational cost}
\label{sec:comp_cost} 
Here, we discuss the computational costs of the GMRES method with different preconditioners. The most computationally intensive operations in the GMRES are matrix-vector products and preconditioning.  Each step of the GMRES thus requires $c_\text{mvp} + c_\text{prec}$, where
\begin{align*}
c_{\text{mvp}}&\text{: cost of matrix-vector products described in eqs.~(\ref{eq:jac-A-mult})--(\ref{eq:jac-B-mult}),}\\
c_{\text{prec}} &\text{: cost of preconditioning.}
\end{align*}
Then the total computational cost of the GMRES is $n_\text{iter}(c_\text{mvp} + c_\text{prec})$, 
where $n_\text{iter}$ refers to the total iteration count. 
The cost of matrix-vector products is largely due to evaluating the first term, 
$\sum_{\ell =1}^{n_{a}} A_{\ell} \delta \bar U H_{\ell }^T$, in  \eqref{eq:jac-A-mult} and, thus, 
the cost can be approximately measured as  $c_{\text{mvp}} \approx n_a (c_{x} + c_{\xi})$, 
where $c_x$ and $c_{\xi}$ are the costs for matrix-matrix products associated 
with $A_\ell$ and $H_\ell$ in the expression $A_{\ell} \delta \bar U H_{\ell }^T$. 
For the preconditioning, we compare two most efficient preconditioners, cMB and truncated chGS with $p_t < p$. 
Let us denote the computational cost of a solve with~$\mathcal M_1$ in~\eqref{eq:MB-matricized-N} by~$c_{\mathcal M_1}$. 
The cMB preconditioner (Algorithm~\ref{alg:NMB}) requires $c_{\text{prec}} =  c_{\mathcal M_1}$ 
and the computational cost of the GMRES with the cMB preconditioner can be approximated~as  
\begin{equation*}
c_{\text{cMB}} = n_{\text{iter}}^{\text{\tiny cMB}} ( \underbrace{n_a (c_x + c_{\xi})}_{c_\text{mvp}} + \underbrace{c_{\mathcal M_1}}_{c_\text{prec}}).
\end{equation*}
The chGS preconditioner (Algorithms \ref{alg:chGS}--\ref{alg:chGS_cont}) requires two truncated matrix-vector products \eqref{eq:matvec-block-N1}--\eqref{eq:matvec-block-N2}, where the truncation is specified by the set $\mathcal I_t$, and applications of the cMB preconditioners for $2p$ times (in the forward and the backward sweep of the Algorithms \ref{alg:chGS}--\ref{alg:chGS_cont}) and, thus, the cost can be assessed as $c_{\text{prec}} \approx 2 n_t (c_x + c_{\xi}) + 2 p\: c_{\mathcal M_1}$, where $n_t = \text{dim}(\mathcal I_t)$. Now we can write the total computational cost of the GMRES method with the chGS preconditioner as
\begin{equation*}
c_{\text{chGS}} = n_{\text{iter}}^{\text{\tiny chGS}} ( \underbrace{n_a(c_x  + c_{\xi})}_{c_\text{mvp}} +  \underbrace{2 n_t (c_x + c_{\xi}) + 2 p\: c_{\mathcal M_1}}_{c_\text{prec}}).
\end{equation*}
From the analytic expressions of the costs, we can see 
that $c_\text{prec}$ for chGS is larger than $c_\text{prec}$ for cMB as chGS 
requires two truncated matrix-vector products $2 n_t (c_x + c_{\xi})$ at each GMRES iteration. 
On the other hand, typically, $n_{\text{iter}}^{\text{cMB}} \gg n_{\text{iter}}^{\text{chGS}}$ and, thus, 
the cMB preconditioner requires more iterations. 
Specifically, the cMB preconditioner needs to perform 
extra $n_{\text{iter}}^{\text{cMB}} - n_{\text{iter}}^{\text{chGS}}$ matrix-vector products, 
with cost $n_a(c_x + c_{\xi})$. To compare the computational costs of the two methods
cMB and chGS($p_t=2$) in practice, we tabulate the values of $n_\xi$, $n_a$ and $n_t$ for varying $m_{\xi}=\{3,5,7\}$ 
and $p=\{3,4,5\}$, see Table~\ref{tab:nant_tbl}.
For problems with coefficients characterized by linear expansion in $\xi$ such as~\eqref{eq:mindlin_exp}, 
cMB could be less expensive since $n_a$ is typically smaller than $n_t$. 
For problems with coefficients characterized by more general (nonlinear) expansions such as~\eqref{eq:a-stoch_exp}, 
chGS with truncated matrix-vector products become more cost efficient because $n_a$ grows exponentially as~$m_{\xi}$ 
and $p$ become larger, whereas $n_t$ remains small. Note that an analogous comparison can be made for chGS and NMB.
%As $m_{\xi}$ and $p$ become larger, however, $n_a$ grows exponentially whereas $n_t$ grows much slower, which makes chGS with truncated matrix-vector products more attractive for high-dimensional problems. 

\begin{table}[!h]
\caption{The number of terms, $n_a$, in the expansion \eqref{eq:stoch-exp-A} modeling linear 
and nonlinear coefficient expansions such as \eqref{eq:a-stoch_exp} and ~\eqref{eq:mindlin_exp}, respectively, 
and the number of terms $n_t$ in the truncation set~$\mathcal I_t$ with $p_t=2$ for varying number of random variables $m_{\xi}$ and the maximum polynomial degree~$p$ of the solution expansion~\eqref{eq:sol_mat}.}
\label{tab:nant_tbl}
\begin{center}
{\footnotesize \renewcommand{\arraystretch}{1.3} 
\begin{tabular}{|c |c |c c c| ccc |ccc |}
\hline
\multicolumn{2}{|c|}{$m_{\xi}$} & \multicolumn{3}{c|}{3} & \multicolumn{3}{c|}{5}& \multicolumn{3}{c|}{7}\\
\hline
\multicolumn{2}{|c|}{$p$} & 3 & 4 & 5 & 3 & 4 & 5& 3 & 4 & 5\\
\hline
\multicolumn{2}{|c|}{$n_{\xi}$} & 20 & 35 & 56 & 56 & 126 & 252 &  120 & 330 & 792\\ 
\hline
\multicolumn{2}{|l|}{$n_a$ (nonlinear)} & 84 & 165 & 286 & 462 & 1287 & 3003& 1716 & 6435 & 19448\\ 
\cline{2-11}
\multicolumn{2}{|l|}{$n_a$ (linear)} &  \multicolumn{3}{c|}{4} & \multicolumn{3}{c|}{6}  & \multicolumn{3}{c|}{8}\\
\hline
\multicolumn{2}{|c|}{$n_t$} & \multicolumn{3}{c|}{10} & \multicolumn{3}{c|}{21}  & \multicolumn{3}{c|}{36}\\
\hline
\end{tabular}
}
\end{center}
\end{table}

\bibliographystyle{siam}
\bibliography{eig_klbs}

\end{document}